\documentclass[a4paper, fleqn]{cas-dc}

\usepackage{amsthm}

\usepackage{subcaption} 
\usepackage{graphicx}
\usepackage{svg}
\usepackage{epstopdf}
\usepackage{caption}
\usepackage{float}
\usepackage{algorithm}
\usepackage{algorithmicx}
\usepackage{algpseudocode}
\usepackage{longtable}
\usepackage[figuresright]{rotating}
\usepackage{pdflscape}
\usepackage[hypcap = true, hypcapspace = 2cm]{caption}
\usepackage{makecell}
\usepackage{multicol}
\usepackage{enumitem}
\usepackage{bigstrut}
\usepackage{tabularx}
\usepackage{threeparttable} 

\usepackage{booktabs}
\usepackage{array, caption, threeparttable}
\newcolumntype{C}[1]{>{\centering\arraybackslash}p{#1}}

\usepackage{amsmath,bm}

\allowdisplaybreaks[4]
\makeatletter
\newenvironment{breakablealgorithm}
{
		\begin{center}
			\refstepcounter{algorithm}
			\hrule height.8pt depth0pt \kern2pt
			\renewcommand{\caption}[2][\relax]{
				{\raggedright\textbf{\ALG@name~\thealgorithm} ##2\par}%
				\ifx\relax##1\relax 
				\addcontentsline{loa}{algorithm}{\protect\numberline{\thealgorithm}##2}%
				\else 
				\addcontentsline{loa}{algorithm}{\protect\numberline{\thealgorithm}##1}%
				\fi
				\kern2pt\hrule\kern2pt
			}
		}{
		\kern2pt\hrule\relax
	\end{center}
}
\makeatother

\usepackage[section]{placeins} 

\usepackage[many]{tcolorbox}
\usepackage{framed}
\usepackage{hyperref} 
\hypersetup{colorlinks=true,
            linkcolor=blue,   
            anchorcolor=black, 
            citecolor=blue,    
            urlcolor=blue      
            }

\usepackage[numbers, sort&compress]{natbib}
\usepackage[numbers]{natbib}

\usepackage{appendix}
\renewcommand\appendix{\par
    \setcounter{section}{0}
    \setcounter{subsection}{0}
    \setlength\LTleft{0pt}
    \setlength\LTright{0pt}
    \captionsetup[table]{margin=-0.03cm}
    \gdef\thesection{Appendix \Alph{section}}}

\def\tsc#1{\csdef{#1}{\textsc{\lowercase{#1}}\xspace}}
\tsc{WGM}
\tsc{QE}

\begin{document}
\let\WriteBookmarks\relax
\renewcommand{\floatpagefraction}{0.8}
\renewcommand{\textfraction}{0.05}
\let\printorcid\relax 
\shorttitle{}



\title[mode = title]{A two-stage stochastic programming framework for oil and gas exploration well portfolio optimization under geological and economic uncertainty}  




\author[1,2,3]{Junyi Cui}[style=chinese]
\author[1, 2]{Junyi Cui}[style=chinese]
\author[3]{Junyi Cui}[style=chinese]
\cormark[1]
\author[1]{Junyi Cui}[style=chinese]
\author[1]{Junyi Cui}[style=chinese]
\author[4]{Junyi Cui}[style=chinese]

\address[1]{School of Sciences, Southwest Petroleum University, Chengdu 610500, China}
\address[2]{Institute for Artificial Intelligence, Southwest Petroleum University, Chengdu 610500, China}
\address[3]{Jagiellonian University in Krakow, Faculty of Mathematics and Computer Science, 30348 Krakow, Poland}

\cortext[1]{Jagiellonian University in Krakow, Faculty of Mathematics and Computer Science, 30348 Krakow, Poland. E-mail address: stanislaw.migorski@uj.edu.pl.} 

\begin{abstract}
Annual oil and gas exploration planning involves selecting a limited portfolio of drilling and appraisal-related projects before geological outcomes are known. This decision is affected by uncertainties in geological success, reserve size, and economic value, while also subject to budget, well-count, success-rate, and reserve-reliability requirements. A strategy based only on expected value is therefore insufficient, as early drilling results may change the value of subsequent follow-up opportunities. This study develops a posterior-informed two-stage stochastic multi-objective optimization framework for exploration well selection under uncertainty. The first stage selects a here-and-now portfolio of frontier traps, appraisal projects, and mature appraisal units. After first-stage outcomes are observed, the second stage determines scenario-dependent recourse projects, including follow-up appraisal, reserve upgrading, conversion-to-proved reserves, rolling extension, and data re-evaluation projects. Geological learning is modeled using a logit-scale posterior updating mechanism that links first-stage success or failure to the success probabilities of related recourse projects. The model maximizes expected net present value and minimizes conditional value-at-risk, while imposing chance constraints on drilling success rate and individual and joint reserve targets. To solve the model, sample average approximation is combined with NSGA-II for first-stage portfolio search and a scenario-wise constrained 0–1 optimization procedure for second-stage evaluation. A numerical case study shows that the proposed framework provides an interpretable risk-return frontier and supports adaptive exploration planning under geological learning, downside-risk control, and reserve-reliability requirements. 
\end{abstract}


\maketitle 

\section{Introduction}\label{sec1}
Oil and gas exploration planning is a high-stakes portfolio decision problem under substantial geological, reserve, and economic uncertainty. In a typical exploration workflow, geoscientists first generate and assess drilling prospects, engineers and economists then translate subsurface assessments into production and financial forecasts, and decision makers finally select a drilling portfolio under limited exploration budgets \cite{bib1}. Risk analysis has therefore been recognized as a central component of petroleum exploration, appraisal, field development, production forecasting, decision making, portfolio management, and real-options analysis \cite{bib2}. Probability of geological success uncertainty is particularly important because the probability of geological success is widely used in prospect risking, expected monetary value analysis, prospect ranking, and exploration portfolio management, but its estimation may be affected by expert judgment and inconsistent geological interpretation \cite{bib3}. Reserve uncertainty is also intrinsic to exploration decisions; the petroleum resources management system represents recoverable resources through uncertainty ranges, including low, best, and high estimates such as P90, P50, and P10. Economic uncertainty is then transmitted through reserve realization, drilling cost, failure loss, NPV, and downside exposure; this is why petroleum decision studies increasingly distinguish expected value from downside risk and upside potential \cite{bib4, bib5}.

Existing studies have provided important foundations for petroleum project selection and portfolio management. Decision analysis \cite{bib6, bib7}, simulation \cite{bib8}, portfolio management \cite{bib9}, and real-options analysis \cite{bib10, bib11, bib12} have been applied in upstream petroleum firms to improve capital allocation and project selection. Exploration portfolio studies have compared rank-and-cut, efficient-frontier, and portfolio-filtering approaches, showing that portfolios based only on risked volume, economic probability of success, or expected monetary value can be useful but may be suboptimal relative to portfolio-level optimization \cite{bib1}. Multi-criteria decision models have also been proposed for selecting oil and gas exploration portfolios, allowing decision makers’ preferences to be represented in a structured way \cite{bib13}. Multi-attribute utility models further recognize that E\&P project selection is stochastic and multi-objective, and that project synergies may affect portfolio value \cite{bib14}. Beyond exploration portfolio selection, stochastic programming has been used in oil and gas field infrastructure planning to handle reserve and production uncertainties that are gradually revealed through investment and operating decisions \cite{bib15, bib16}. These studies collectively demonstrate that petroleum investment planning is uncertain, constrained, and multi-objective.

However, several limitations remain. First, deterministic and ranking-based approaches often summarize uncertainty through expected or risked metrics, while portfolio-level reserve reliability is not explicitly controlled; LaCosta et al. \cite{bib13} showed that rank-and-cut portfolios can perform well but remain suboptimal compared with optimized efficient-frontier portfolios. Second, project-by-project evaluation may overlook portfolio diversification and resource-allocation effects, which are central to portfolio management under capital constraints \cite{bib4}. Third, multi-criteria and utility-based models can represent preferences and multiple attributes, but they are usually formulated as single-stage portfolio selection models rather than adaptive exploration strategies \cite{bib13, bib14}. Fourth, sequential exploration studies have shown that drilling outcomes may provide valuable geological information and reshape subsequent drilling policies. For example, Martinelli et al. \citep{bib17} combine Bayesian networks with dynamic programming to derive exploration strategies in which later drilling decisions depend on earlier well outcomes, while Martinelli et al. \citep{bib18} extend this idea to graphical models and approximate dynamic programming. These studies primarily focus on sequential prospect selection and information updating, whereas annual exploration portfolio planning requires the simultaneous treatment of budget limits, well-count requirements, reserve chance constraints, downside-risk control, and category-specific second-stage recourse decisions \cite{bib19, bib20}. Finally, second-stage follow-up actions are often simplified as local ranking or myopic decisions, whereas constrained recourse actions in annual planning are combinatorial because follow-up appraisal, reserve upgrading, and conversion-to-proved-reserve projects compete for the same remaining budget and drilling workload.

To address these gaps, this paper develops a posterior-informed two-stage stochastic programming framework for annual oil and gas exploration well selection. The first stage determines a here-and-now exploration portfolio before geological outcomes are observed. After first-stage outcomes are realized, the second stage selects scenario-dependent recourse projects, including follow-up appraisal wells, conversion-to-proved-reserve projects, reserve upgrading projects, rolling extension projects, and data re-evaluation projects. The posterior updating mechanism is formulated on the logit scale. This choice is motivated by the Bayesian odds interpretation: posterior odds can be written as prior odds multiplied by an evidence factor, and the logarithmic transformation converts evidence accumulation into an additive update \citep{bib21}. In this study, first-stage success or failure provides geological evidence that shifts the log-odds of related second-stage projects, while the logistic transformation maps the updated value back to a valid probability. This provides a parsimonious way to connect geological learning with recourse decisions without allowing probabilities to move outside the admissible range.

The proposed model optimizes two conflicting objectives: maximizing expected NPV and minimizing the conditional value-at-risk of downside losses. CVaR is adopted because it focuses on tail losses rather than total variability, which is consistent with the managerial concern that exploration portfolios may generate severe losses under unfavorable geological outcomes. The theoretical basis of coherent risk measurement was established by Artzner et al. \citep{bib22}, and CVaR was subsequently formulated as an optimization-compatible tail-risk measure by \citep{bib23}. CVaR has also been widely used in portfolio optimization with objective functions and constraints \citep{bib24}. In the petroleum context, downside-oriented risk measures are particularly relevant because decision makers may care more about losses below a benchmark than about symmetric variability around the mean \citep{bib5}. The proposed formulation therefore combines expected economic value with downside-risk control, while also imposing chance constraints on drilling success rate, individual reserve targets, and joint reserve target satisfaction.

A hybrid scenario-based solution framework is designed to solve the resulting model. Sample average approximation is used to evaluate expected NPV, empirical chance constraints, and CVaR under finite geological-success, reserve, and economic-value scenarios, following the established framework for stochastic discrete optimization \citep{bib25}. NSGA-II is used to search the first-stage binary portfolio space and approximate the risk-return Pareto frontier, taking advantage of its non-dominated sorting and crowding-distance mechanisms for multi-objective optimization \citep{bib26}. For each first-stage portfolio and each realized geological scenario, the second-stage decision is evaluated through an explicit 0--1 recourse subproblem solved by a non-dominated-state dynamic programming procedure. This framework separates strategic first-stage portfolio search from scenario-wise recourse optimization, while preserving the adaptive structure of two-stage decision making. In other words, second-stage actions are evaluated as constrained recourse portfolios rather than as independent follow-up projects selected only by local ranking.

The main contributions of this study are fourfold:
\begin{itemize}
	\item It develops an uncertainty characterization scheme for exploration portfolio optimization that links geological success probability, reserve-size uncertainty, and economic outcomes within a unified stochastic decision framework.
	
	\item It proposes a posterior-informed two-stage stochastic programming model in which first-stage drilling outcomes update second-stage project probabilities and activate scenario-dependent recourse decisions.
	
	\item It incorporates risk-aware and chance-constrained portfolio requirements by jointly considering expected NPV, CVaR, drilling success reliability, individual reserve targets, and joint reserve target satisfaction.
	
	\item It designs a hybrid solution framework that combines NSGA-II-based first-stage portfolio search with scenario-wise 0--1 recourse evaluation, enabling the generation of interpretable Pareto-efficient exploration strategies under operational and reserve-management constraints.
\end{itemize}

The remainder of this paper is organized as follows. Section \ref{sec2} characterizes the geological success, reserve, and economic uncertainties considered in the exploration portfolio problem. Section \ref{sec3} formulates the posterior-informed two-stage stochastic programming model. Section \ref{sec4} presents the scenario generation, SAA, NSGA-II search procedure, and exact recourse evaluation method. Section \ref{sec5} reports the case study and computational results. Section \ref{sec6} discusses managerial implications and model extensions.

\section{Uncertainty characterization}\label{sec2}
Oil and gas exploration planning is made before drilling outcomes are fully observed. The uncertainty considered in this study is therefore characterized at the project level and then propagated to the two-stage portfolio model. Specifically, three stochastic components are considered: geological success uncertainty, reserve uncertainty, and economic-value uncertainty. This section defines these stochastic primitives, while their decision-dependent use in the optimization model and their finite-scenario implementation are presented in Sections~\ref{sec3} and~\ref{sec4}, respectively.

\subsection{Geological success uncertainty}
Geological success uncertainty describes whether a candidate project can achieve a technically successful drilling outcome. For a generic project \(q\), the prior success probability is denoted by \(p_q^0\). When the success probability is assessed through geological risk factors, a multiplicative prospect-risking representation is used:
\begin{flalign}\label{1}
	& p_q^0
	=
	P_q^{\mathrm{src}}
	P_q^{\mathrm{res}}
	P_q^{\mathrm{trap}}
	P_q^{\mathrm{pre}}
	P_q^{\mathrm{mig}}, &
\end{flalign}
where \(P_q^{\mathrm{src}}\), \(P_q^{\mathrm{res}}\), \(P_q^{\mathrm{trap}}\), \(P_q^{\mathrm{pre}}\), and \(P_q^{\mathrm{mig}}\) denote the probabilities associated with source rock effectiveness, reservoir development, trap validity, hydrocarbon preservation, and migration--charge effectiveness, respectively. 

For the first-stage project \(i\in I\), let \(\xi_i^s=1\) if project \(i\) succeeds in scenario \(s\), and \(\xi_i^s=0\) otherwise. Then
\begin{flalign}\label{2}
	& \xi_i^s\in\{0,1\},
	\;
	\mathbb{P}(\xi_i^s=1)=p_i^0,
	\; i\in I,\;s\in S. &
\end{flalign}

Unlike first-stage projects, recourse projects are evaluated after the first-stage drilling outcomes have been observed. Their success probabilities are therefore updated on the logit scale. This treatment is consistent with the odds form of Bayesian updating, in which new evidence modifies prior odds into posterior odds, and taking logarithms converts multiplicative evidence effects into additive log-odds terms \cite{bib27}. It is also related to weights-of-evidence modelling in geoscience, where multiple evidence layers are combined through additive evidence weights \cite{bib28}. In petroleum exploration, Bayesian updating has been used to integrate geological chance of success with geophysical evidence and obtain revised success probabilities \cite{bib29}. Following this logic, the observed first-stage outcomes are represented as a net log-odds evidence term for each recourse project.

Let \(\mathcal{A}_j\) denote the set of first-stage projects that provide geological information for recourse project \(j\in J\). The net log-odds evidence received by project \(j\) in scenario \(s\) is
\begin{flalign}\label{3}
	& \Delta_j^s(\boldsymbol{x})
	=
	\sum_{i\in\mathcal{A}_j}
	\theta_{ij}x_i(2\xi_i^s-1), &
\end{flalign}
where \(\theta_{ij}\) is the logit-scale information strength from first-stage project \(i\) to recourse project \(j\). A successful first-stage outcome increases the log-odds of related recourse projects when \(\theta_{ij}>0\), whereas a failed outcome decreases it. If project \(i\) is not selected in the first stage, then \(x_i=0\), and its outcome does not contribute to the posterior evidence.

The posterior success probability of recourse project \(j\) in scenario \(s\) is then expressed as
\begin{flalign}\label{4}
	& p_j^s(\boldsymbol{x})
	=
	\mathcal{P}_{[\underline{p},\overline{p}]}
	\left[
	\sigma\left(
	\mathrm{logit}(p_j^0)
	+
	\Delta_j^s(\boldsymbol{x})
	\right)
	\right], &
\end{flalign}
where \(\sigma(a)=1/(1+\exp(-a))\), \(\mathrm{logit}(p)=\ln[p/(1-p)]\), and \(\mathcal{P}_{[\underline{p},\overline{p}]}[\cdot]\) projects the updated probability onto an admissible interval. Conditional on this posterior probability, the second-stage success indicator satisfies
\begin{flalign}\label{5}
	& \mathbb{P}
	\left(
	\zeta_j^{s,k}=1
	\mid \boldsymbol{x},\boldsymbol{\xi}^s
	\right)
	=
	p_j^s(\boldsymbol{x}),
	\quad j\in J. &
\end{flalign}

\subsection{Reserve uncertainty}
Reserve uncertainty reflects the uncertainty in geological scale and reservoir quality conditional on project success. For each project, oil and gas reserve potentials are determined by uncertain volumetric parameters, including area, effective thickness, porosity, water saturation, and formation volume factor. These inputs are represented by three-point estimates corresponding to the lower, median, and upper quantile values, denoted here as \((P10,P50,P90)\) for notational convenience.

For project \(q\), the random oil reserve potential is calculated by the volumetric method:
\begin{flalign}\label{6}
	& \widetilde{R}_q^o
	=
	\kappa_o
	\frac{
		\widetilde{A}_q^o
		\widetilde{h}_q^o
		\widetilde{\phi}_q^o
		(1-\widetilde{S}_{w,q}^o)
	}{
		\widetilde{B}_{o,q}
	}, &
\end{flalign}
and the random gas reserve potential is
\begin{flalign}\label{7}
	& \widetilde{R}_q^g
	=
	\kappa_g
	\frac{
		\widetilde{A}_q^g
		\widetilde{h}_q^g
		\widetilde{\phi}_q^g
		(1-\widetilde{S}_{w,q}^g)
	}{
		\widetilde{B}_{g,q}
	}, &
\end{flalign}
where \(A\), \(h\), \(\phi\), \(S_w\), and \(B\) denote area, effective thickness, porosity, water saturation, and formation volume factor, respectively; \(\kappa_o\) and \(\kappa_g\) are unit conversion coefficients.

The reserve contribution of project \(q\) to reserve indicator \(m\) is obtained by applying classification coefficients to the sampled oil and gas reserve potentials:
\begin{flalign}\label{8}
	& \widetilde{r}_{qm}
	=
	\lambda_{qm}^o \widetilde{R}_q^o
	+
	\lambda_{qm}^g \widetilde{R}_q^g. &
\end{flalign}

The coefficients \(\lambda_{qm}^o\) and \(\lambda_{qm}^g\) allow the same physical reserve potential to contribute differently to predicted, controlled, and proved reserve indicators. Thus, reserve uncertainty enters the optimization model through both economic value and chance-constrained reserve target satisfaction.
\begin{figure*}[htbp]
	\centering
	\includegraphics[width=\linewidth]{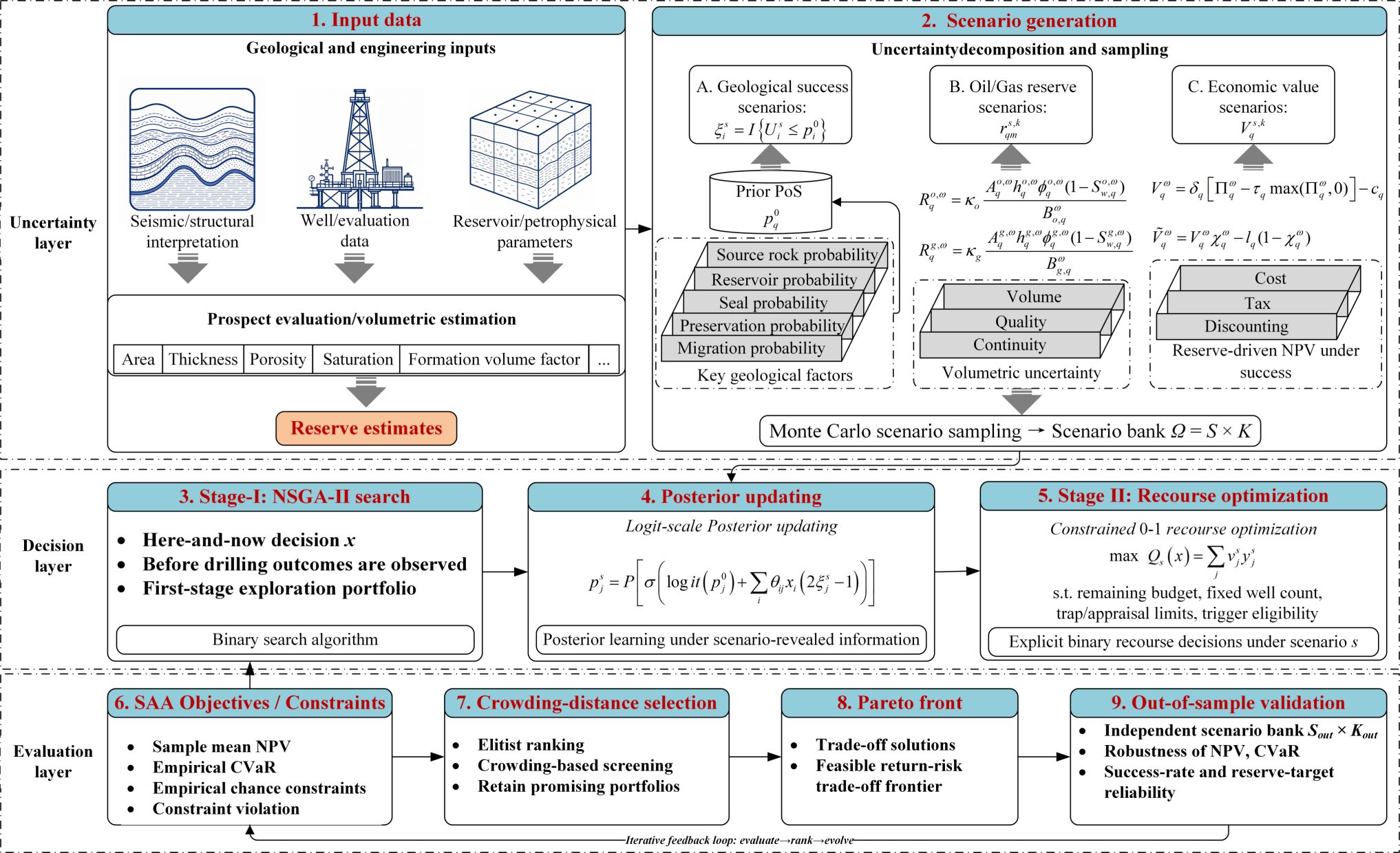} 
	\caption{Framework of the proposed two-stage stochastic multi-objective optimization method.}
	\label{fig1}
\end{figure*}

\subsection{Economic-value uncertainty}
In this study, economic-value uncertainty is induced by geological success and reserve realization. Price, unit production cost, fixed operating cost, tax rate, discount factor, and economic reserve coefficients are treated as project-level planning parameters, while the realized economic value remains stochastic because the success state and reserve potential are uncertain.

Let \(a_q^\omega\) denote the success indicator of project \(q\) under realization \(\omega\), where \(a_q^\omega=\xi_i^s\) for \(q=i\in I\) and \(a_q^\omega=\zeta_j^{s,k}\) for \(q=j\in J\). The project-level realized payoff is
\begin{flalign}\label{9}
	& \mathcal{V}_q^\omega
	=
	V_q^\omega a_q^\omega
	-
	l_q(1-a_q^\omega), &
\end{flalign}
where \(V_q^\omega\) denotes the success-state NPV specified in Section~\ref{sec3}, and \(l_q\) is the net loss in the failure state.

This representation distinguishes high-upside projects from low-risk projects in a direct way. A project with large reserve potential may create high economic value when successful, but it can still generate a loss if drilling fails. Conversely, an appraisal-related project may have lower upside but more stable reserve and economic contributions. The resulting trade-off between expected value and downside exposure is captured in the bi-objective two-stage model developed in Section~\ref{sec3}.

The three uncertainty components above enter the optimization model through success indicators, posterior success probabilities, reserve contributions, and realized project values. Their finite-scenario realization is described in Section~\ref{sec4}.

\section{Two-stage stochastic programming model}\label{sec3}

This section formulates the annual exploration well-selection problem as a two-stage stochastic programming model. The first-stage decision determines the here-and-now exploration portfolio before drilling outcomes are observed. After first-stage outcomes are realized, the second-stage decision selects scenario-dependent recourse projects, including follow-up appraisal wells, conversion-to-proved-reserve projects, reserve upgrading projects, rolling extension projects, and data re-evaluation projects. The uncertainty characterization in Section~\ref{sec2} enters the model through success indicators, posterior success probabilities, reserve contributions, and realized economic values. The framework is illustrated in Fig.~\ref{fig1}.
\begin{table}[t]
	\centering
	\footnotesize
	\setlength{\tabcolsep}{0pt}
	\renewcommand{\arraystretch}{1.10}
	\caption{Sets, indices, decisions, and uncertainty variables}
	\label{Tab1}
	\begin{tabularx}{\columnwidth}{@{}>{\raggedright\arraybackslash}p{0.28\columnwidth}
			>{\raggedright\arraybackslash}X@{}}
		\hline
		Symbol & Description \\
		\hline
		\(I\) & Set of first-stage candidate projects. \\
		\(J\) & Set of second-stage candidate recourse projects. \\
		\(q\) & Generic project index, \(q\in I\cup J\). \\
		\(s\) & First-stage scenario index. \\
		\(k\) & Conditional second-stage sub-scenario index. \\
		\(S\) & Set of first-stage scenarios. \\
		\(K\) & Set of second-stage sub-scenarios. \\
		\(\Omega\) & Set of scenario pairs, \(\Omega=\{(s,k)\mid s\in S,k\in K\}\). \\
		\(M\) & Set of reserve indicators. \\
		\(\omega\) & Generic realization index; \(\omega=s\) for first-stage projects and \(\omega=(s,k)\) for second-stage projects. \\
		\(M^{\mathrm{act}}\) & Set of reserve indicators with activated chance constraints. \\
		\(I^{\mathrm{fix}}\) & Set of mandatory first-stage projects. \\
		\(I^{\mathrm{trap}},J^{\mathrm{trap}}\) & Sets of trap-related projects in the first and second stages. \\
		\(I^{\mathrm{app}},J^{\mathrm{app}}\) & Sets of appraisal-related projects in the first and second stages. \\
		\(a_q^\omega\) & Generic success indicator of project \(q\) under realization \(\omega\). \\
		\(\mathcal{A}_j\) & Set of first-stage projects providing geological information for recourse project \(j\). \\
		\(\mathcal{A}_j^{+},\mathcal{A}_j^{-},\mathcal{A}_j^{0}\) & Success-, failure-, and unconditional-trigger sets for recourse project \(j\). \\
		\(x_i\) & Binary first-stage selection decision for project \(i\). \\
		\(y_j^s\) & Binary second-stage recourse decision for project \(j\) in scenario \(s\). \\
		\(\xi_i^s\) & Success indicator of first-stage project \(i\) in scenario \(s\). \\
		\(\zeta_j^{s,k}\) & Success indicator of second-stage project \(j\) under \((s,k)\). \\
		\(e_j^s(\boldsymbol{x})\) & Eligibility indicator of recourse project \(j\) in scenario \(s\). \\
		\(p_q^0\) & Prior success probability of project \(q\). \\
		\(p_j^s(\boldsymbol{x})\) & Posterior success probability of recourse project \(j\) in scenario \(s\). \\
		\(\theta_{ij}\) & Logit-scale information strength from project \(i\) to project \(j\). \\
		\(\Delta_j^s(\boldsymbol{x})\) & Net log-odds evidence received by recourse project \(j\). \\
		\hline
	\end{tabularx}
\end{table}

\begin{table}[t]
	\centering
	\footnotesize
	\setlength{\tabcolsep}{0pt}
	\renewcommand{\arraystretch}{1.10}
	\caption{Economic, constraint, and recourse-related parameters}
	\label{Tab2}
	\begin{tabularx}{\columnwidth}{@{}>{\raggedright\arraybackslash}p{0.28\columnwidth}
			>{\raggedright\arraybackslash}X@{}}
		\hline
		Symbol & Description \\
		\hline
		\(\pi_q^o,\pi_q^g\) & Oil and gas prices of project \(q\). \\
		\(u_q^o,u_q^g\) & Unit production costs of oil and gas for project \(q\). \\
		\(\eta_q^o,\eta_q^g\) & Economic reserve coefficients for oil and gas. \\
		\(f_q\) & Fixed operating cost of project \(q\). \\
		\(\tau_q\) & Income tax rate of project \(q\). \\
		\(\delta_q\) & Discount factor of project \(q\). \\
		\(\lambda_{qm}^o,\lambda_{qm}^g\) & Classification coefficients mapping oil and gas reserve potentials to reserve indicator \(m\). \\
		\(\kappa_o,\kappa_g\) & Unit conversion coefficients for oil and gas volumetric reserve calculation. \\
		\(\underline{p},\overline{p}\) & Lower and upper bounds of admissible posterior success probabilities. \\
		\(c_q\) & Investment cost of project \(q\). \\
		\(w_q\) & Number of wells associated with project \(q\). \\
		\(l_q\) & Failure loss of project \(q\). \\
		\(R_q^{o,\omega},R_q^{g,\omega}\) & Oil and gas reserve potentials of project \(q\) under realization \(\omega\). \\
		\(r_{qm}^{\omega}\) & Reserve contribution of project \(q\) to indicator \(m\). \\
		\(V_q^{\omega}\) & Success-state NPV of project \(q\) under realization \(\omega\). \\
		\(Z^{s,k}(\boldsymbol{x},\boldsymbol{y}^s)\) & Realized annual NPV under scenario pair \((s,k)\). \\
		\(L^{s,k}(\boldsymbol{x},\boldsymbol{y}^s)\) & Downside loss under scenario pair \((s,k)\). \\
		\(\mathcal{R}_m^{s,k}(\boldsymbol{x},\boldsymbol{y}^s)\) & Total reserve contribution to indicator \(m\) under \((s,k)\). \\
		\(\Gamma^{s,k}(\boldsymbol{x},\boldsymbol{y}^s)\) & Scenario-dependent drilling success rate. \\
		\(B_1,N_1\) & First-stage investment and well-count limits. \\
		\(B,N\) & Annual total investment limit and annual well-count requirement. \\
		\(B^{\mathrm{trap}},B^{\mathrm{app}}\) & Investment limits for trap-related and appraisal-related projects. \\
		\(H_m\) & Minimum target for reserve indicator \(m\). \\
		\(\alpha_m\) & Required satisfaction probability for reserve indicator \(m\). \\
		\(\alpha^{\mathrm{joint}}\) & Required joint satisfaction probability for active reserve targets. \\
		\(\rho^{\min}\) & Minimum required drilling success rate. \\
		\(\alpha^{\mathrm{sr}}\) & Required satisfaction probability for the success-rate constraint. \\
		\(\beta\) & Confidence level for CVaR calculation. \\
		\(C_1(\boldsymbol{x})\) & First-stage investment cost. \\
		\(W_1(\boldsymbol{x})\) & First-stage well count. \\
		\(\bar{B}(\boldsymbol{x})\) & Remaining annual investment capacity. \\
		\(\bar{N}(\boldsymbol{x})\) & Remaining annual well-count requirement. \\
		\(\nu_j^s(\boldsymbol{x})\) & Expected economic value of recourse project \(j\). \\
		\(v_j^s(\boldsymbol{x})\) & Recourse value including reserve-shortfall compensation. \\
		\(\gamma\) & Weight of reserve-shortfall compensation. \\
		\(Q_s(\boldsymbol{x})\) & Optimal value of the scenario-wise recourse subproblem. \\
		\hline
	\end{tabularx}
\end{table}

\subsection{Sets, indices and parameters}\label{sec:notation}
Let \(I\) denote the set of first-stage candidate projects, including frontier traps, appraisal projects, and mature appraisal units. Let \(J\) denote the set of second-stage candidate recourse projects. The first-stage geological scenario is indexed by \(s\in S\), and the conditional second-stage sub-scenario is indexed by \(k\in K\). The set of reserve indicators is defined as
\begin{flalign}\label{10}
	& M=\{\mathrm{po},\mathrm{pg},\mathrm{co},\mathrm{cg},
	\mathrm{ro},\mathrm{rg}\}, &
\end{flalign}
where \(\mathrm{po}\), \(\mathrm{pg}\), \(\mathrm{co}\), \(\mathrm{cg}\), \(\mathrm{ro}\), and \(\mathrm{rg}\) denote predicted oil reserves, predicted gas reserves, controlled oil reserves, controlled gas reserves, proved oil reserves, and proved gas reserves, respectively. Let \(M^{\mathrm{act}}\subseteq M\) denote the set of reserve indicators whose chance constraints are activated in the annual planning problem.

For compactness, define the successful pre-tax operating profit of project \(q\in I\cup J\) under realization \(\omega\) as
\begin{flalign}\label{11}
	&\begin{aligned}
		\Pi_q^\omega
		={}&
		(\pi_q^o-u_q^o)\eta_q^o R_q^{o,\omega}
		+(\pi_q^g-u_q^g)\eta_q^g R_q^{g,\omega}
		-f_q .
	\end{aligned}&
\end{flalign}

The successful NPV is then
\begin{flalign}\label{12}
	& V_q^\omega
	=
	\delta_q
	\left[
	\Pi_q^\omega-\tau_q\max(\Pi_q^\omega,0)
	\right]
	-c_q , &
\end{flalign}
where \(R_q^{o,\omega}\) and \(R_q^{g,\omega}\) are the sampled oil and gas reserve potentials; \(\pi_q^o\) and \(\pi_q^g\) are oil and gas prices; \(u_q^o\) and \(u_q^g\) are unit production costs; \(\eta_q^o\) and \(\eta_q^g\) are economic reserve coefficients; \(f_q\) is the fixed operating cost; \(\tau_q\) is the income tax rate; \(\delta_q\) is the discount factor; and \(c_q\) is the investment cost. If project \(q\) fails, its realized value is \(-l_q\), where \(l_q\) is the failure loss. The main notation used in the model is summarized in Table~\ref{Tab1} and \ref{Tab2}.

\subsection{Two-stage stochastic programming formulation}
\label{sec:model}

The proposed model determines an annual exploration portfolio under geological, reserve, and economic uncertainty. The first-stage decision is made before drilling outcomes are observed, whereas the second-stage decision is implemented after the first-stage geological information has been revealed. Thus, the model explicitly distinguishes the here-and-now portfolio decision from the scenario-dependent recourse decision.

For each first-stage project \(i\in I\), the binary decision variable is defined as
\begin{flalign}\label{13}
	& x_i=
	\begin{cases}
		1, & \text{if project } i \text{ is selected in the first stage},\\
		0, & \text{otherwise}.
	\end{cases} &
\end{flalign}

For each second-stage project \(j\in J\) and first-stage scenario \(s\in S\), the recourse decision is
\begin{flalign}\label{14}
	& y_j^s=
	\begin{cases}
		1, & \text{if project } j \text{ is selected as recourse in scenario } s,\\
		0, & \text{otherwise}.
	\end{cases} &
\end{flalign}

After the first-stage outcomes are observed, a second-stage project can be selected only when its trigger condition is satisfied. Let \(\mathcal{A}_j^{+}\), \(\mathcal{A}_j^{-}\), and \(\mathcal{A}_j^{0}\) denote the first-stage projects that trigger recourse project \(j\) through success, failure, and unconditional linkage, respectively. The eligibility indicator is
\begin{flalign}\label{15}
	&\begin{aligned}
		e_j^s(\boldsymbol{x})
		=\mathbb{I}\Bigg\{
		\sum_{i\in\mathcal{A}_j^{+}}x_i\xi_i^s
		+\sum_{i\in\mathcal{A}_j^{-}}x_i(1-\xi_i^s) 
		+\sum_{i\in\mathcal{A}_j^{0}}x_i
		\ge 1
		\Bigg\}.
	\end{aligned}&
\end{flalign}

The posterior success probability of recourse project \(j\) in scenario \(s\) is denoted by \(p_j^s(\boldsymbol{x})\), which is obtained from the logit-scale posterior updating rule described in Section~\ref{sec2}. Conditional on \(p_j^s(\boldsymbol{x})\), the second-stage success indicator satisfies
\begin{flalign}\label{16}
	& \mathbb{P}
	\left(
	\zeta_j^{s,k}=1
	\mid \boldsymbol{x},\boldsymbol{\xi}^s
	\right)
	=
	p_j^s(\boldsymbol{x}),
	\qquad j\in J. &
\end{flalign}

The realized annual NPV under scenario pair \((s,k)\) is
\begin{flalign}\label{17}
	&\begin{aligned}
		Z^{s,k}(\boldsymbol{x},\boldsymbol{y}^s)
		={}&
		\sum_{i\in I}
		x_i\left[
		V_i^s\xi_i^s-l_i(1-\xi_i^s)
		\right] \\
		&+
		\sum_{j\in J}
		y_j^s\left[
		V_j^{s,k}\zeta_j^{s,k}
		-l_j(1-\zeta_j^{s,k})
		\right].
	\end{aligned}&
\end{flalign}
Here, \(V_i^s\) and \(V_j^{s,k}\) are the success-state NPVs determined by the sampled reserve and economic parameters, while \(l_i\) and \(l_j\) are failure losses.

The total reserve contribution to indicator \(m\in M\) is
\begin{flalign}\label{18}
	&\begin{aligned}
		\mathcal{R}_m^{s,k}(\boldsymbol{x},\boldsymbol{y}^s)
		={}&
		\sum_{i\in I}
		x_i\xi_i^s r_{im}^s
		+
		\sum_{j\in J}
		y_j^s\zeta_j^{s,k} r_{jm}^{s,k}.
	\end{aligned}&
\end{flalign}

The scenario-dependent drilling success rate is defined as
\begin{flalign}\label{19}
	& \Gamma^{s,k}(\boldsymbol{x},\boldsymbol{y}^s)
	=
	\frac{
		\sum_{i\in I}w_i x_i\xi_i^s
		+
		\sum_{j\in J}w_j y_j^s\zeta_j^{s,k}
	}{
		\sum_{i\in I}w_i x_i
		+
		\sum_{j\in J}w_j y_j^s
	}. &
\end{flalign}

The first objective maximizes the expected annual NPV:
\begin{flalign}\label{20}
	& \max\; F_1
	=
	\frac{1}{|S||K|}
	\sum_{s\in S}\sum_{k\in K}
	Z^{s,k}(\boldsymbol{x},\boldsymbol{y}^s). &
\end{flalign}

To control downside economic exposure, the loss variable is defined relative to the zero-profit benchmark:
\begin{flalign}\label{21}
	& L^{s,k}(\boldsymbol{x},\boldsymbol{y}^s)
	=
	\max\left\{
	-Z^{s,k}(\boldsymbol{x},\boldsymbol{y}^s),\;0
	\right\}. &
\end{flalign}

For confidence level \(\beta\in(0,1)\), the CVaR of downside losses is defined by the Rockafellar-Uryasev \citep{bib23, bib30} representation:
\begin{flalign}\label{22}
	& \mathrm{CVaR}_{\beta}(L)
	=
	\min_{\alpha\in\mathbb{R}}
	\left\{
	\alpha+
	\frac{1}{1-\beta}
	\mathbb{E}\left[(L-\alpha)^+\right]
	\right\}, &
\end{flalign}
where \(\alpha\) is the value-at-risk threshold and \((a)^+=\max\{a,0\}\). Under the finite scenario representation, the second objective is written as
\begin{flalign}\label{23}
	&\min \;{F_2} = \mathop {\min }\limits_{_{\alpha  \in \mathbb{R}}} \left\{ {\alpha  + \frac{{\sum\limits_{s \in S} {\sum\limits_{k \in K} {{{\left[ {{L^{s,k}}({\boldsymbol{x}},{{\boldsymbol{y}}^s}) - \alpha } \right]}^ + }} } }}{{(1 - \beta )|S||K|}}} \right\}.
	&
\end{flalign}

The first-stage investment and well-count constraints are
\begin{flalign}\label{24}
	& \sum_{i\in I}c_i x_i \le B_1, &
\end{flalign}
and
\begin{flalign}\label{25}
	& \sum_{i\in I}w_i x_i \le N_1. &
\end{flalign}

For each scenario \(s\), the annual total investment and well-count requirements are
\begin{flalign}\label{26}
	& \sum_{i\in I}c_i x_i
	+
	\sum_{j\in J}c_j y_j^s
	\le B,
	\; s\in S, &
\end{flalign}
and
\begin{flalign}\label{27}
	& \sum_{i\in I}w_i x_i
	+
	\sum_{j\in J}w_j y_j^s
	=
	N,
	\; s\in S. &
\end{flalign}

The investment limits for trap-related projects and appraisal projects are
\begin{flalign}\label{28}
	& \sum_{i\in I^{\mathrm{trap}}}c_i x_i
	+
	\sum_{j\in J^{\mathrm{trap}}}c_j y_j^s
	\le B^{\mathrm{trap}},
	\; s\in S, &
\end{flalign}
and
\begin{flalign}\label{29}
	& \sum_{i\in I^{\mathrm{app}}}c_i x_i
	+
	\sum_{j\in J^{\mathrm{app}}}c_j y_j^s
	\le B^{\mathrm{app}},
	\; s\in S. &
\end{flalign}

The recourse decision must satisfy the trigger eligibility condition:
\begin{flalign}\label{30}
	& y_j^s \le e_j^s(\boldsymbol{x}),
	\; j\in J,\; s\in S. &
\end{flalign}

The drilling success-rate chance constraint is
\begin{flalign}\label{31}
	& \frac{1}{|S||K|}
	\sum_{s\in S}\sum_{k\in K}
	\mathbb{I}
	\left\{
	\Gamma^{s,k}(\boldsymbol{x},\boldsymbol{y}^s)
	\ge \rho^{\min}
	\right\}
	\ge
	\alpha^{\mathrm{sr}}. &
\end{flalign}

Let \(M^{\mathrm{act}}\subseteq M\) denote the set of active reserve targets. For each reserve indicator \(m\in M^{\mathrm{act}}\), the individual reserve chance constraint is
\begin{flalign}\label{32}
	& \frac{1}{|S||K|}
	\sum_{s\in S}\sum_{k\in K}
	\mathbb{I}
	\left\{
	\mathcal{R}_m^{s,k}(\boldsymbol{x},\boldsymbol{y}^s)
	\ge H_m
	\right\}
	\ge
	\alpha_m,
	\; m\in M^{\mathrm{act}}. &
\end{flalign}

The joint reserve chance constraint requires all active reserve targets to be satisfied simultaneously:
\begin{flalign}\label{33}
	& \frac{1}{|S||K|}
	\sum_{s\in S}\sum_{k\in K}
	\mathbb{I}
	\left\{
	\mathcal{R}_m^{s,k}(\boldsymbol{x},\boldsymbol{y}^s)
	\ge H_m,\;
	\forall m\in M^{\mathrm{act}}
	\right\}
	\ge
	\alpha^{\mathrm{joint}}. &
\end{flalign}

If a subset of first-stage projects \(I^{\mathrm{fix}}\subseteq I\) must be included for operational or contractual reasons, the following fixed-selection constraint is imposed:
\begin{flalign}\label{34}
	& x_i=1,
	\; i\in I^{\mathrm{fix}}. &
\end{flalign}

Finally, all first-stage and second-stage decisions are binary:
\begin{flalign}\label{35}
	& x_i\in\{0,1\},
	\; i\in I,
	\;
	y_j^s\in\{0,1\},
	\; j\in J,\;s\in S. &
\end{flalign}

Equations~\eqref{20}--\eqref{35} define a two-stage stochastic bi-objective portfolio optimization model. The first objective evaluates the expected economic value of the annual exploration plan, while the second objective controls the tail risk of negative NPV outcomes. The recourse variables \(\boldsymbol{y}^s\) are not fixed in advance; they are optimized after scenario-specific geological information is revealed, subject to remaining budget, well-count, project-category, trigger-eligibility, and reserve-reliability requirements.

\section{Methodology}\label{sec4}
The model in Section~\ref{sec3} contains binary first-stage decisions, scenario-dependent recourse decisions, posterior probability updating, chance constraints, and a CVaR-based risk objective. Directly solving the full stochastic model becomes computationally demanding as the number of candidate projects and uncertainty scenarios increases. Therefore, this study adopts a simulation-based solution framework that combines SAA, NSGA-II, and scenario-wise explicit recourse evaluation.

The framework separates the strategic and adaptive components of the problem. NSGA-II searches over first-stage portfolios, while each candidate portfolio is evaluated by posterior updating and a constrained 0--1 recourse subproblem in every first-stage scenario. Thus, the second-stage decision is not treated as a fixed add-on or a local greedy ranking, but as a scenario-dependent recourse portfolio subject to remaining investment, well-count, and project-category constraints.

\subsection{Scenario generation}\label{sec:scenario_generation}
A finite scenario bank is generated to represent the stochastic components defined in Section~\ref{sec2}. For first-stage project \(i\in I\), the success indicator in scenario \(s\) is generated as
\begin{flalign}\label{36}
	& \xi_i^s=\mathbb{I}\{U_i^s\le p_i^0\},
	\; i\in I,\ s\in S, &
\end{flalign}
where \(U_i^s\sim U(0,1)\). For recourse project \(j\in J\), the second-stage success indicator is generated after the posterior probability \(p_j^s(\boldsymbol{x})\) has been updated:
\begin{flalign}\label{37}
	& \zeta_j^{s,k}(\boldsymbol{x})
	=
	\mathbb{I}\{U_j^{s,k}\le p_j^s(\boldsymbol{x})\},
	\; j\in J,\ s\in S,\ k\in K, &
\end{flalign}

For notational simplicity, the model formulation uses \(\zeta_j^{s,k}\) to denote this posterior-dependent realization once \(\boldsymbol{x}\) and \(\boldsymbol{\xi}^s\) are fixed. The random numbers \(U_i^s\) and \(U_j^{s,k}\) are generated once and reused for all candidate portfolios. This common-random-number design ensures that different first-stage portfolios are compared under the same sampled uncertainty, thereby reducing simulation noise in Pareto ranking.

For a generic uncertain volumetric parameter \(\psi_q\), the sampled value is written as
\begin{flalign}\label{38}
	&\psi_q^\omega
	\sim
	\mathcal{D}
	\left(
	\psi_q^{0.10},
	\psi_q^{0.50},
	\psi_q^{0.90}
	\right), &
\end{flalign}
where \(\omega=s\) for first-stage projects and \(\omega=(s,k)\) for second-stage projects. Positive-valued volumetric parameters are sampled using a lognormal distribution calibrated by the three quantiles when feasible. Bounded parameters, such as porosity and water saturation, are clipped to physically admissible intervals. If lognormal calibration is not appropriate, a PERT-type distribution is used. The sampled volumetric parameters are then transformed into \(R_q^{o,\omega}\), \(R_q^{g,\omega}\), \(r_{qm}^{\omega}\), and \(V_q^{\omega}\) according to the definitions in Sections~\ref{sec2} and~\ref{sec3}.

\subsection{Sample average approximation}\label{sec:saa}
The stochastic objective functions and chance constraints are evaluated over a finite in-sample scenario bank. Let
\begin{flalign}\label{39}
	& \Omega
	=
	\{(s,k)\mid s\in S,\; k\in K\},
	\qquad
	N_{\Omega}=|S||K|. &
\end{flalign}
For a given first-stage portfolio \(\boldsymbol{x}\), the explicit recourse evaluation procedure determines an optimal recourse decision \(\boldsymbol{y}^{s*}\) for each first-stage scenario \(s\in S\). The expected annual NPV is then approximated by
\begin{flalign}\label{40}
	& \widehat{F}_1(\boldsymbol{x})
	=
	\frac{1}{N_{\Omega}}
	\sum_{(s,k)\in\Omega}
	Z^{s,k}(\boldsymbol{x},\boldsymbol{y}^{s*}). &
\end{flalign}

For a scenario-wise event \(\mathcal{E}^{s,k}\), its empirical probability is estimated as
\begin{flalign}\label{41}
	& \widehat{\mathbb{P}}(\mathcal{E})
	=
	\frac{1}{N_{\Omega}}
	\sum_{(s,k)\in\Omega}
	\mathbb{I}\{\mathcal{E}^{s,k}\}. &
\end{flalign}
Accordingly, the drilling success-rate constraint, individual reserve chance constraints, and joint reserve chance constraint are evaluated by empirical frequencies over \(\Omega\).

Let \(L^{s,k}\) be the downside loss defined in Eq.~\eqref{21}. In the SAA implementation, the CVaR objective is evaluated by the finite-sample counterpart of the Rockafellar-Uryasev representation:
\begin{flalign}\label{42}
	&{\widehat{{\rm{CVaR}}}_\beta }({\boldsymbol{x}}) = \mathop {\min }\limits_{\alpha  \in \mathbb{R}} \left\{ {\alpha  + \frac{{\sum\limits_{(s,k) \in \Omega } {{{\left[ {{L^{s,k}}({\boldsymbol{x}},{{\boldsymbol{y}}^{s*}}) - \alpha } \right]}^ + }} }}{{(1 - \beta ){N_\Omega }}}} \right\},&
\end{flalign}

This expression is consistent with the CVaR definition in Eq.~\eqref{22} and keeps the risk objective in the same finite-scenario form as the expected-value objective. The bi-objective evaluation used by NSGA-II is written in minimization form as
\begin{flalign}\label{43}
	& \boldsymbol{F}^{\mathrm{SAA}}(\boldsymbol{x})
	=
	\left(
	-\widehat{F}_1(\boldsymbol{x}),
	\widehat{\mathrm{CVaR}}_{\beta}(\boldsymbol{x})
	\right). &
\end{flalign}

\subsection{NSGA-II for first-stage portfolio search}\label{sec:nsga}
Each candidate solution in NSGA-II is encoded as a binary first-stage portfolio:
\begin{flalign}\label{44}
	& \boldsymbol{x}
	=
	(x_1,x_2,\ldots,x_{|I|}),
	\; x_i\in\{0,1\}. &
\end{flalign}

For each chromosome, mandatory projects are repaired by enforcing the fixed-selection requirement in Eq.~\eqref{34}. The candidate portfolio is then evaluated over the scenario bank. In each first-stage scenario \(s\), the algorithm observes \(\xi_i^s\), updates \(e_j^s(\boldsymbol{x})\) and \(p_j^s(\boldsymbol{x})\), solves the recourse subproblem to obtain \(\boldsymbol{y}^{s*}\), and finally computes the sample-based objectives and constraint violations.

The objective vector used by NSGA-II is
\begin{flalign}\label{45}
	& \boldsymbol{F}(\boldsymbol{x})
	=
	\left(
	-\widehat{F}_1(\boldsymbol{x}),
	\widehat{\mathrm{CVaR}}_{\beta}(\boldsymbol{x})
	\right). &
\end{flalign}

Constraint handling follows a feasibility-first rule. Let \(g_h(\boldsymbol{x})\le 0\) denote the \(h\)-th deterministic or empirical chance constraint, and let \(\mathcal{H}\) be the index set of all such constraints.
\begin{flalign}\label{46}
	& \mathrm{Viol}(\boldsymbol{x})
	=
	\sum_{h\in\mathcal{H}}
	\max\{g_h(\boldsymbol{x}),0\}. &
\end{flalign}

Feasible solutions are ranked by non-dominated sorting and crowding distance. Infeasible solutions are compared according to their aggregate violation. The output of NSGA-II is a set of feasible non-dominated first-stage portfolios and their associated scenario-dependent recourse policies.

\subsection{Explicit recourse evaluation}\label{sec:exact_recourse}

For a fixed first-stage portfolio \(\boldsymbol{x}\) and first-stage scenario \(s\), the eligible recourse set is
\begin{flalign}\label{47}
	& J^s(\boldsymbol{x})
	=
	\{j\in J\mid e_j^s(\boldsymbol{x})=1\}, &
\end{flalign}
the remaining annual capacities are
\begin{flalign}\label{48}
	& \bar{B}(\boldsymbol{x})=B-C_1(\boldsymbol{x}),
	\;
	\bar{N}(\boldsymbol{x})=N-W_1(\boldsymbol{x}). &
\end{flalign}

The remaining trap-related and appraisal-related investment capacities are
\begin{flalign}\label{49}
	&\begin{aligned}
		\bar{B}^{\mathrm{trap}}(\boldsymbol{x})
		&=
		B^{\mathrm{trap}}
		-\sum_{i\in I^{\mathrm{trap}}}c_ix_i,\\
		\bar{B}^{\mathrm{app}}(\boldsymbol{x})
		&=
		B^{\mathrm{app}}
		-\sum_{i\in I^{\mathrm{app}}}c_ix_i .
	\end{aligned}&
\end{flalign}

Let the first-stage reserve contribution already obtained in scenario \(s\) be
\begin{flalign}\label{50}
	& R_{1m}^{s}(\boldsymbol{x})
	=
	\sum_{i\in I}r_{im}^{s}\xi_i^sx_i,
	\; m\in M^{\mathrm{act}}, &
\end{flalign}
the expected economic value of recourse project \(j\) is
\begin{flalign}\label{51}
	& \nu_j^s(\boldsymbol{x})
	=
	p_j^s(\boldsymbol{x})
	\frac{1}{|K|}\sum_{k\in K}V_j^{s,k}
	-
	[1-p_j^s(\boldsymbol{x})]l_j . &
\end{flalign}

To guide recourse selection toward unmet reserve targets, a reserve-shortfall term is added. The recourse value is defined as
\begin{flalign}\label{52}
	&\begin{aligned}
		v_j^s(\boldsymbol{x})
		={}&
		\nu_j^s(\boldsymbol{x})
		+\gamma c_j
		\sum_{m\in M^{\mathrm{act}}}
		\mathbb{I}\{R_{1m}^{s}(\boldsymbol{x})<H_m\} \\
		&\times
		\frac{
			p_j^s(\boldsymbol{x})
			\left(|K|^{-1}\sum_{k\in K}r_{jm}^{s,k}\right)
		}{H_m},
	\end{aligned}&
\end{flalign}
where \(\gamma\) controls the weight of reserve-shortfall compensation.

The second-stage recourse portfolio is obtained by solving the following 0--1 subproblem:
\begin{flalign}\label{53}
	& Q_s(\boldsymbol{x})
	=
	\max_{\boldsymbol{y}^s}
	\sum_{j\in J^s(\boldsymbol{x})}
	v_j^s(\boldsymbol{x})y_j^s, &
\end{flalign}
subject to
\begin{flalign}\label{54}
	& \sum_{j\in J^s(\boldsymbol{x})}c_jy_j^s
	\le \bar{B}(\boldsymbol{x}), &
\end{flalign}
\begin{flalign}\label{55}
	& \sum_{j\in J^s(\boldsymbol{x})}w_jy_j^s
	= \bar{N}(\boldsymbol{x}), &
\end{flalign}
\begin{flalign}\label{56}
	& \sum_{j\in J^{\mathrm{trap}}\cap J^s(\boldsymbol{x})}
	c_jy_j^s
	\le \bar{B}^{\mathrm{trap}}(\boldsymbol{x}), &
\end{flalign}
\begin{flalign}\label{57}
	& \sum_{j\in J^{\mathrm{app}}\cap J^s(\boldsymbol{x})}
	c_jy_j^s
	\le \bar{B}^{\mathrm{app}}(\boldsymbol{x}), &
\end{flalign}
and
\begin{flalign}\label{58}
	& y_j^s\in\{0,1\},
	\; j\in J^s(\boldsymbol{x}). &
\end{flalign}

This subproblem is a constrained binary recourse portfolio problem. It is solved by a non-dominated-state dynamic programming procedure. Each state records accumulated total investment, trap-related investment, appraisal-related investment, well count, recourse value, and selected project indices. For states with the same well count, state \(a\) dominates state \(b\) if
\begin{flalign}\label{59}
	& C_a\le C_b,\;
	C_a^{\mathrm{trap}}\le C_b^{\mathrm{trap}},\;
	C_a^{\mathrm{app}}\le C_b^{\mathrm{app}},\;
	Q_a\ge Q_b. &
\end{flalign}

Dominated states are removed because they cannot lead to a better feasible recourse portfolio. Let \(\mathcal{D}_t\) be the set of non-dominated states after processing \(t\) eligible recourse projects. The state transition is
\begin{flalign}\label{60}
	& \mathcal{D}_{t+1}
	=
	\mathrm{ND}
	\left(
	\mathcal{D}_{t}
	\cup
	\{d\oplus j_{t+1}:d\in\mathcal{D}_{t}\}
	\right), &
\end{flalign}
where \(d\oplus j_{t+1}\) denotes the feasible state obtained by adding project \(j_{t+1}\), and \(\mathrm{ND}(\cdot)\) denotes non-dominated pruning. After all eligible projects have been processed, the best feasible state gives
\begin{flalign}\label{61}
	& \boldsymbol{y}^{s*}
	\in
	\arg\max_{\boldsymbol{y}^s}
	\left\{
	\sum_{j\in J^s(\boldsymbol{x})}
	v_j^s(\boldsymbol{x})y_j^s
	\; \middle| \;
	\text{constraints } \eqref{54}\text{--}\eqref{58}
	\right\}, &
\end{flalign}
where $ s\in S $, when all non-dominated states are retained, this procedure solves the specified recourse subproblem exactly. Therefore, the adaptive second-stage action is evaluated as a constrained portfolio decision rather than as an independent ranking of recourse projects. Here, exactness refers to the solution of the scenario-wise binary recourse subproblem, rather than to an exact deterministic-equivalent solution of the full stochastic bi-objective model.

The complete computational procedure is summarized in Algorithm~\ref{alg1}. After the feasible non-dominated portfolios are obtained from the in-sample scenario bank, the exported solutions are re-evaluated using an independent out-of-sample scenario bank. This out-of-sample evaluation is used to check the stability of expected NPV, CVaR, drilling success reliability, and reserve target satisfaction.
\begin{breakablealgorithm}
	\caption{Hybrid SAA--NSGA-II framework with scenario-wise recourse evaluation}
	\label{alg1}
	\begin{algorithmic}[1]
		\State Generate in-sample common random numbers, reserve scenarios, and economic-value samples.
		\State Initialize a binary NSGA-II population of first-stage portfolios.
		\For{each generation}
		\For{each candidate portfolio \(\boldsymbol{x}\)}
		\State Repair mandatory first-stage projects.
		\For{each first-stage scenario \(s\in S\)}
		\State Observe first-stage outcomes \(\xi_i^s\).
		\State Update \(e_j^s(\boldsymbol{x})\) and \(p_j^s(\boldsymbol{x})\).
		\State Solve the 0--1 recourse subproblem to obtain \(\boldsymbol{y}^{s*}\).
		\EndFor
		\State Evaluate \(\widehat{F}_1(\boldsymbol{x})\), \(\widehat{\mathrm{CVaR}}_{\beta}(\boldsymbol{x})\), and constraint violations.
		\EndFor
		\State Apply non-dominated sorting, crowding-distance selection, crossover, and mutation.
		\EndFor
		\State Export feasible non-dominated portfolios and their recourse policies.
		\State Re-evaluate exported solutions using independent out-of-sample scenarios.
	\end{algorithmic}
\end{breakablealgorithm}

\section{Experiments and result analysis}\label{sec5}
The preceding sections formulate the proposed two-stage stochastic multi-objective optimization model and describe its solution procedure based on scenario sampling, posterior-informed recourse evaluation, and NSGA-II search. This section examines whether the proposed framework can produce economically meaningful and risk-aware exploration portfolios under realistic uncertainty. The computational experiments are organized to answer five questions: whether the stochastic model generates a clear return--risk trade-off frontier, whether posterior-informed recourse improves decision quality, whether explicit recourse optimization is preferable to rule-based recourse selection, whether the obtained portfolios remain stable under out-of-sample scenarios, and how sensitive the results are to risk preference, geological learning strength, and scenario sampling.

\subsection{Experimental setting and input characterization}
The numerical study is constructed from a synthetic but geologically consistent exploration portfolio. The candidate set contains 30 first-stage projects and 50 second-stage appraisal projects. The first-stage projects represent here-and-now drilling decisions made before exploration outcomes are observed, whereas the second-stage appraisal projects represent recourse opportunities that can be selected after scenario-specific geological information has been revealed. Each project is characterized by prior geological success probability, drilling cost, well-count requirement, reserve-related parameters, and economic value parameters. These inputs are used consistently in scenario generation, posterior updating, recourse optimization, and performance evaluation.

The stochastic experiment uses an in-sample scenario bank for optimization and an independent out-of-sample scenario bank for validation. In the in-sample evaluation, \(S=200\) first-stage scenarios are generated, and each first-stage scenario is associated with \(K=20\) second-stage sub-scenarios. For out-of-sample validation, an independent scenario bank with \(S_{\mathrm{out}}=1000\) and \(K_{\mathrm{out}}=20\) is generated using a different random seed, see Table.~\ref{Tab3} for details. This design separates portfolio search from external validation and reduces the risk that the reported performance is driven by a particular sampled scenario set.
\begin{table}[htbp]
	\caption{Computational setting of the numerical experiments.}\label{Tab3}
	\centering
	\begin{tabular}{ll}
		\hline
		Item & Setting \\
		\hline
		First-stage candidate projects & 30 \\
		Second-stage appraisal projects & 50 \\
		In-sample first-stage scenarios & \(S=200\) \\
		Second-stage sub-scenarios & \(K=20\) \\
		Out-of-sample first-stage scenarios & \(S_{\mathrm{out}}=1000\) \\
		Out-of-sample sub-scenarios & \(K_{\mathrm{out}}=20\) \\
		NSGA-II population size & 100 \\
		NSGA-II generations & 500 \\
		Exported feasible non-dominated portfolios & 27 \\
		\hline
	\end{tabular}
\end{table}

\begin{figure*}[!t]
	\centering	
	\includegraphics[width=0.9\textwidth]{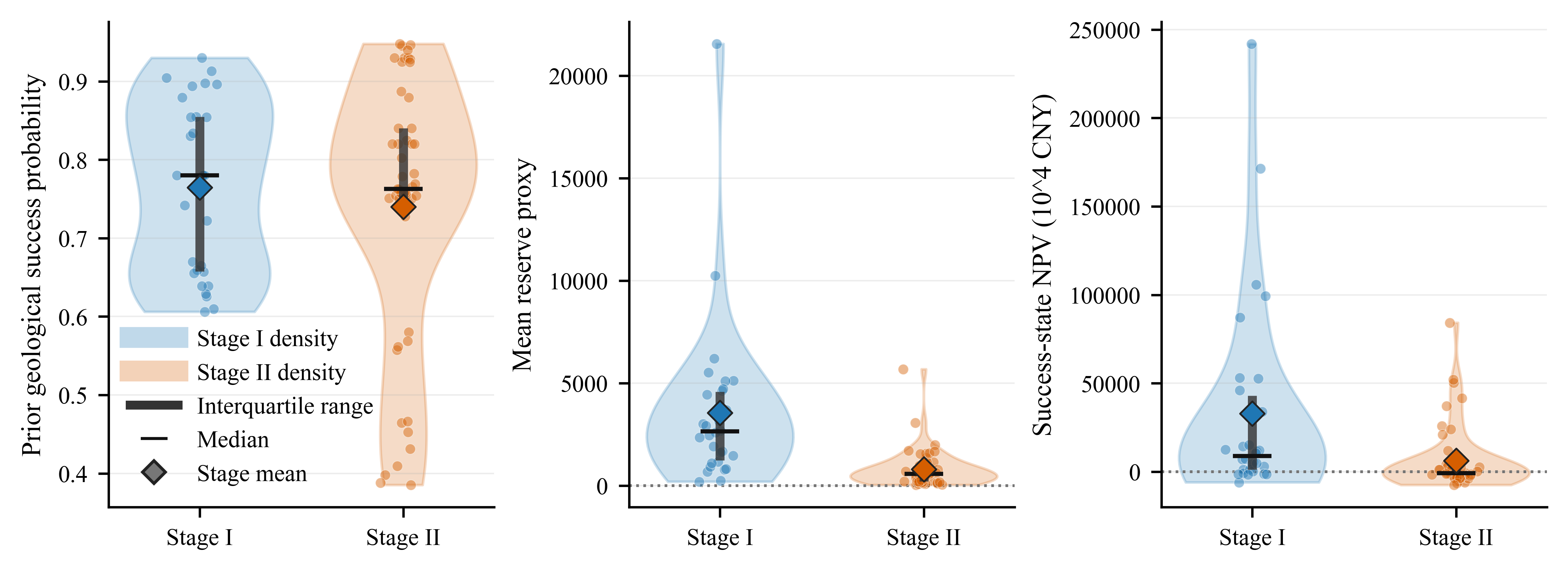}
	\caption{Input distributions of key project-level parameters for first-stage projects and second-stage appraisal projects.}
	\label{fig2}
\end{figure*}

Fig.~\ref{fig2} summarizes the main input distributions for the first- and second-stage candidate sets. The violin density describes the empirical distribution of each parameter, the thick vertical segment indicates the interquartile range, the horizontal marker denotes the median, and the diamond denotes the mean. Compared with second-stage appraisal projects, first-stage projects generally have larger drilling costs, larger reserve potential, and higher success-state NPV dispersion. This structure is consistent with the intended decision hierarchy: the first stage contains larger and more strategic exploration commitments, while the second stage provides more flexible appraisal and compensatory opportunities after geological information is revealed.

The parameter distributions also explain why a two-stage stochastic structure is required. High-value projects are not necessarily low-risk projects, and projects with attractive reserve potential may also require larger drilling commitments. Therefore, portfolio quality cannot be evaluated only by deterministic expected value. In the following subsections, each candidate portfolio is evaluated by its sample-average expected NPV, empirical CVaR-based downside risk, success-probability reliability, reserve-reliability feasibility, and out-of-sample robustness. These indicators jointly describe the economic, risk, and geological reliability dimensions of exploration decision-making.

\subsection{Pareto frontier and deterministic benchmark comparison}
This subsection reports the main optimization results generated by the proposed two-stage stochastic multi-object\-ive model and then compares them with a deterministic mean-value benchmark. The Pareto frontier obtained from the proposed model is the primary computational result, while the deterministic benchmark is used only as a reference to examine the consequence of ignoring uncertainty during portfolio optimization.

\begin{figure*}
	\centering
	\includegraphics[width=0.7\textwidth]{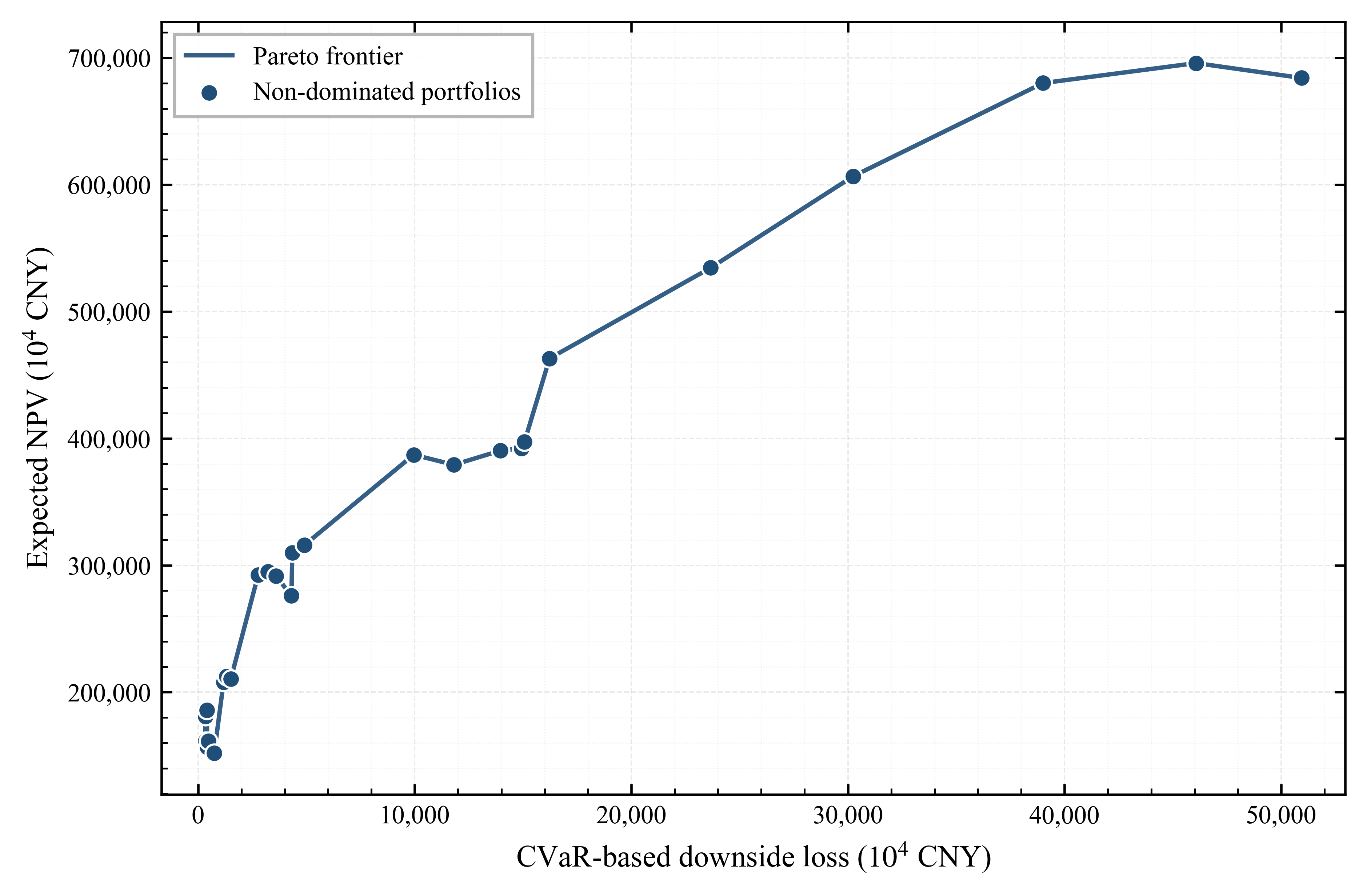}
	\caption{Pareto frontier obtained by the proposed two-stage stochastic multi-objective optimization model.}
	\label{fig3}
\end{figure*}

Figure~\ref{fig3} presents the Pareto frontier directly obtained from the proposed two-stage stochastic model. The algorithm returns 27 feasible non-dominated portfolios. These portfolios are not determined by a pre-specified number of exported solutions; instead, they are retained according to feasibility and dominance relationships. The frontier shows a clear trade-off between expected NPV and CVaR-based downside risk, indicating that higher expected return is generally associated with greater exposure to adverse scenarios.

Table.~\ref{Tab4} summarizes the in-sample and out-of-sample performance of the 27 non-dominated portfolios. In the out-of-sample evaluation, the expected NPV ranges from \(151{,}967.48\) to \(695{,}884.35\) \((10^4\ \mathrm{CNY})\), while the CVaR-based downside loss ranges from \(337.62\) to \(50{,}934.16\) \-\((10^4\ \mathrm{CNY})\). The out-of-sample success reliability remains high, ranging from 0.9994 to 1.0000. By contrast, the joint reserve reliability varies from 0.6060 to 0.9275, which shows that geological success and reserve deliverability are not equivalent decision criteria. Therefore, reserve-related reliability must be explicitly evaluated in exploration portfolio optimization.
\begin{table*}[!t]
	\centering
	\caption{Performance summary of the Pareto portfolios generated by the proposed model.}
	\label{Tab4}
	\begin{tabular}{llcccc}
		\hline
		Sample & Statistic & ENPV & CVaR loss & Success reliability & Reserve reliability \\
		\hline
		In-sample & Minimum & 155356.94 & 0.00 & 0.9972 & 0.6292 \\
		In-sample & Mean    & 346395.45 & 3038.74 & 0.9996 & 0.7742 \\
		In-sample & Maximum & 718457.78 & 19017.25 & 1.0000 & 0.9498 \\
		Out-of-sample & Minimum & 151967.48 & 337.62 & 0.9994 & 0.6060 \\
		Out-of-sample & Mean    & 339920.36 & 11191.13 & 0.9999 & 0.7931 \\
		Out-of-sample & Maximum & 695884.35 & 50934.16 & 1.0000 & 0.9275 \\
		\hline
	\end{tabular}
	\\
	\vspace{1mm}
	\footnotesize{Note: ENPV and CVaR loss are reported in \(10^4\ \mathrm{CNY}\).}
\end{table*}

The difference between in-sample and out-of-sample results is also informative. The mean expected NPV decreases from \(346{,}395.45\) to \(339{,}920.36\) \((10^4\ \mathrm{CNY})\), whereas the mean CVaR loss increases from \(3{,}038.74\) to \(11{,}191.13\) \((10^4\ \mathrm{CNY})\). This pattern suggests that the central economic performance is relatively stable, but downside risk becomes more visible when the portfolios are evaluated under an independent and larger scenario bank. This supports the use of out-of-sample validation for assessing tail-risk robustness.

\begin{figure}
	\centering
	\includegraphics[width=\linewidth]{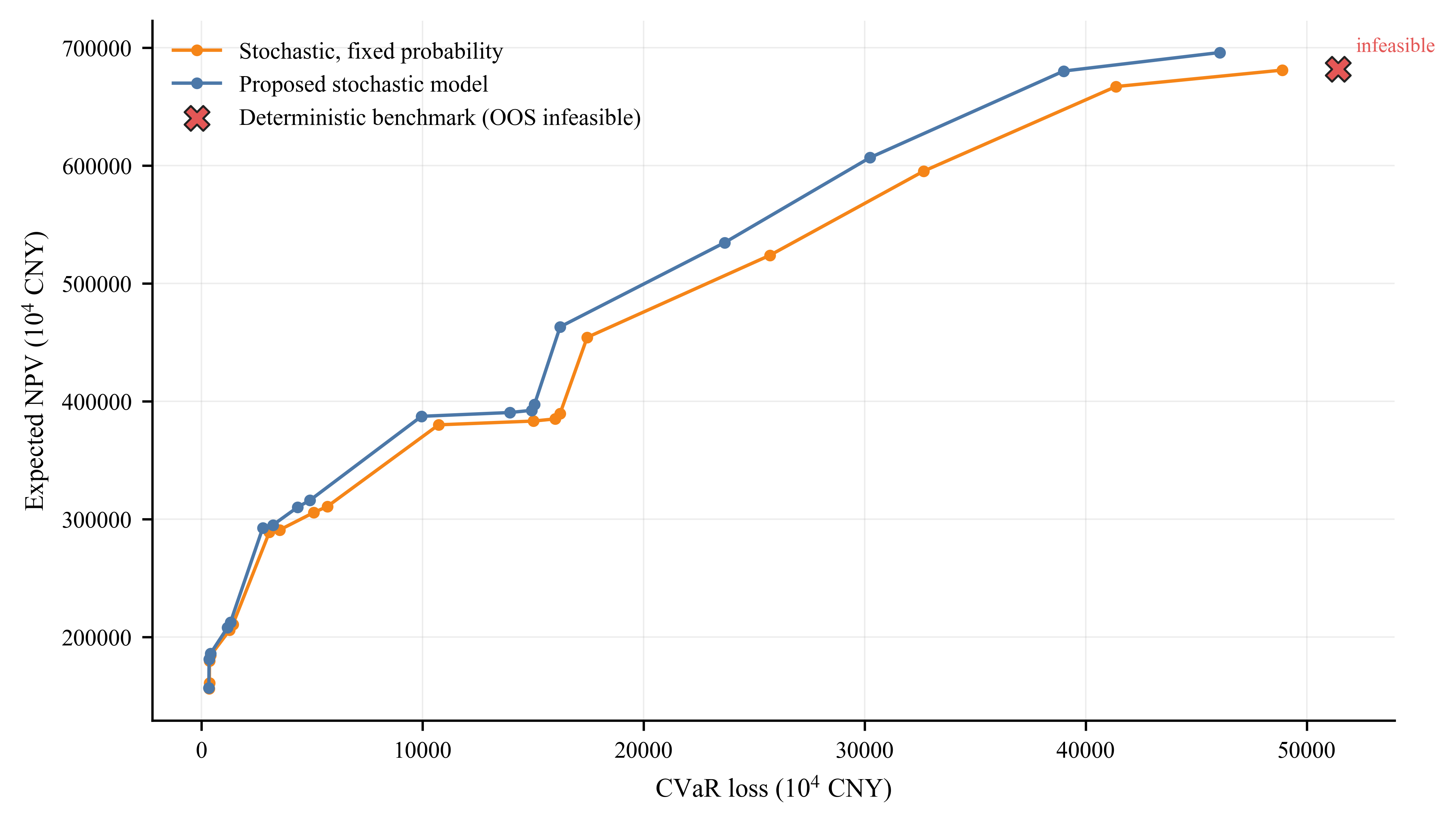}
	\caption{Out-of-sample comparison between the stochastic Pareto frontier and the deterministic mean-value benchmark.}
	\label{fig4}
\end{figure}

To further examine the value of stochastic optimization, a deterministic mean-value benchmark is constructed. The deterministic benchmark optimizes the portfolio using mean parameter values and is then evaluated under the same out-of-sample stochastic scenario bank as the stochastic portfolios. Figure~\ref{fig4} compares the deterministic benchmark with the stochastic Pareto frontier.

The deterministic benchmark obtains a high out-of-sample expected NPV of \(681{,}792.61\) \((10^4\ \mathrm{CNY})\), which is close to the upper end of the stochastic frontier, as shown in Table.~\ref{Tab5}. However, this apparent economic attractiveness is accompanied by a CVaR loss of \(51{,}418.29\) \((10^4\ \mathrm{CNY})\) , a success reliability of 0.9960, and a positive constraint violation. Therefore, the deterministic portfolio is not robust once it is exposed to stochastic geological, reserve, and economic realizations.
\begin{table*}[!t]
	\centering
	\caption{Out-of-sample comparison between stochastic portfolios and the deterministic benchmark.}
	\label{Tab5}
	\begin{tabular}{lccccc}
		\hline
		Portfolio type & ENPV & CVaR loss & Success reliability & Reserve reliability & Feasibility \\
		\hline
		Stochastic Pareto, minimum & 151967.48 & 337.62 & 0.9994 & 0.6060 & feasible \\
		Stochastic Pareto, mean & 339920.36 & 11191.13 & 0.9999 & 0.7931 & feasible \\
		Stochastic Pareto, maximum & 695884.35 & 50934.16 & 1.0000 & 0.9275 & feasible \\
		Deterministic benchmark & 681792.61 & 51418.29 & 0.9960 & 0.8319 & infeasible \\
		\hline
	\end{tabular}
	\\
	\vspace{1mm}
	\footnotesize{Note: ENPV and CVaR loss are reported in \(10^4\ \mathrm{CNY}\).}
\end{table*}

\begin{figure}[!t]
	\centering
	\includegraphics[width=\linewidth]{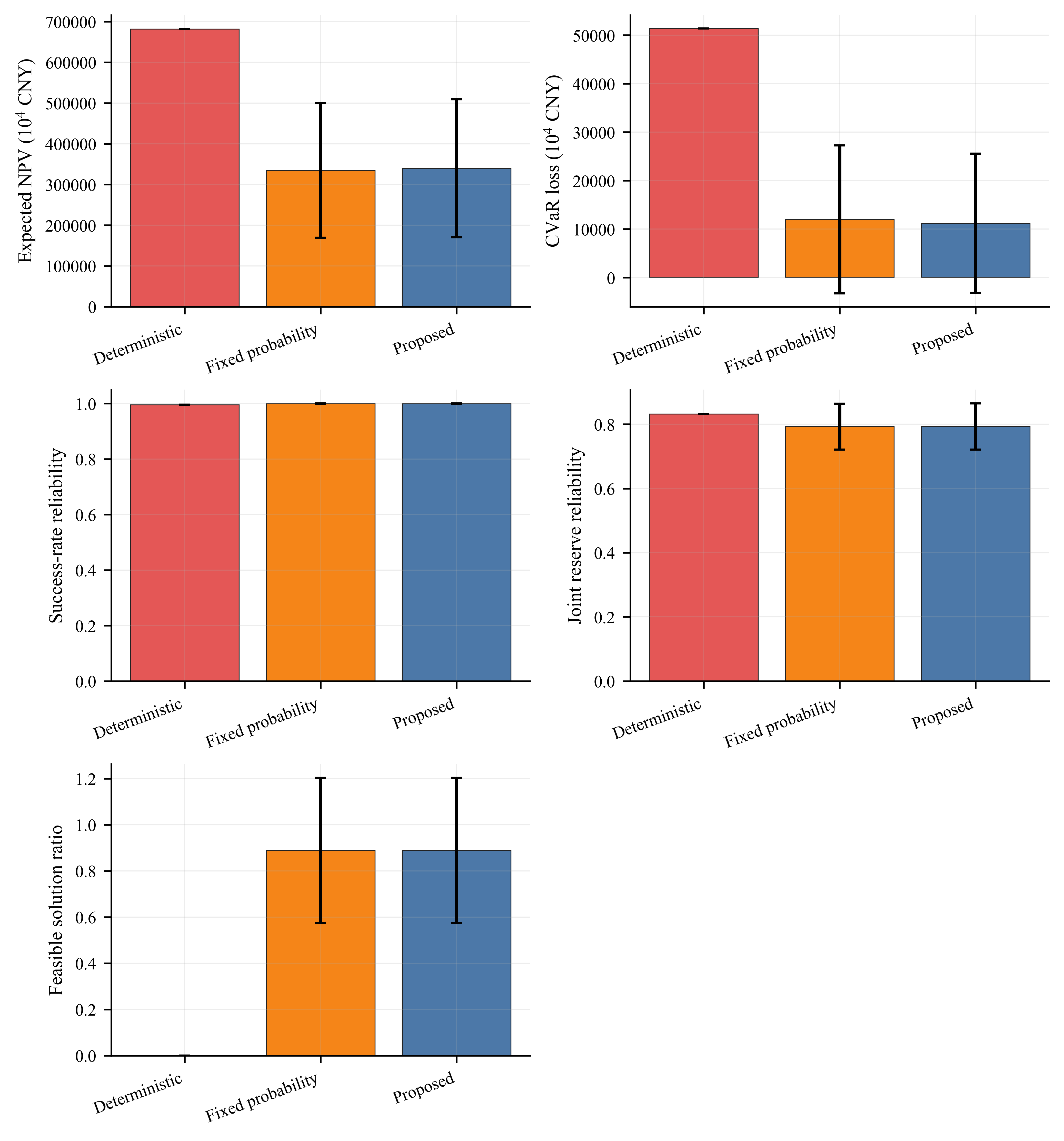}
	\caption{Metric-level comparison between stochastic portfolios and the deterministic benchmark under out-of-sample scenarios.}
	\label{fig5}
\end{figure}

Fig.~\ref{fig5} provides a metric-level comparison. The deterministic benchmark emphasizes nominal expected value but performs poorly in downside-risk and feasibility terms. In contrast, the stochastic portfolios provide a set of feasible alternatives with different return--risk profiles. This result supports the necessity of the proposed stochastic formulation: exploration portfolio optimization should not be based only on deterministic expected value, but should jointly consider expected return, downside risk, geological success reliability, reserve reliability, and recourse feasibility.

\subsection{Value of posterior-informed and explicit recourse optimization}
This subsection evaluates the two key recourse-related components of the proposed framework: posterior-informed probability updating and explicit second-stage recourse optimization. Three recourse settings are first compared: no recourse, exact recourse with fixed prior probabilities, and exact recourse with posterior-updated probabilities. The comparison isolates the value of geological learning from the value of simply allowing second-stage decisions.

\begin{figure}[!t]
	\centering
	\includegraphics[width=\linewidth]{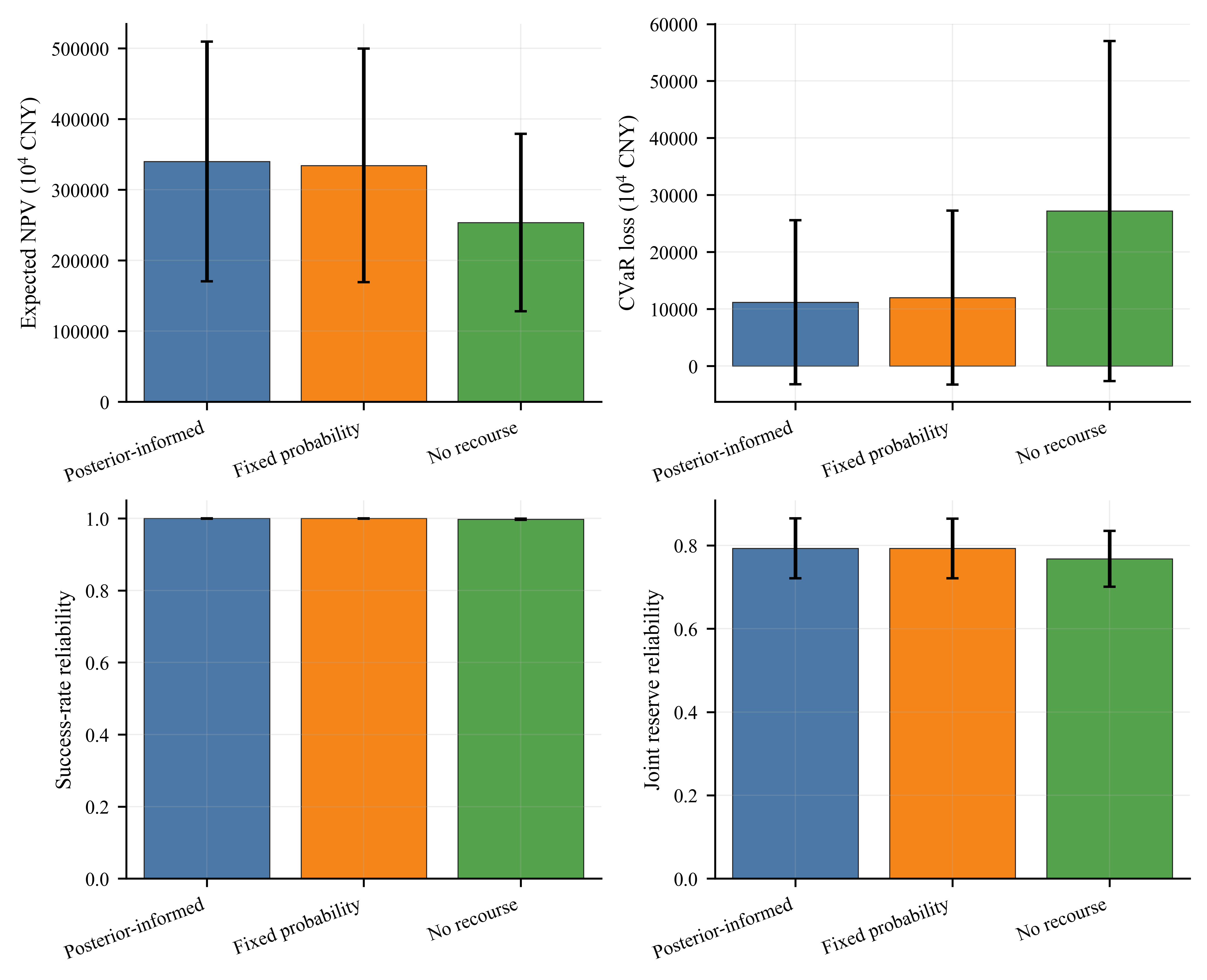}
	\caption{Out-of-sample value of posterior-informed recourse compared with fixed-probability recourse and no recourse.}
	\label{fig6}
\end{figure}

Fig.~\ref{fig6} shows that introducing recourse substantially improves out-of-sample performance. Compared with the no-recourse setting, posterior-informed exact recourse increases the mean expected NPV from \(253{,}732.31\) to \(339{,}920.36\) \((10^4\ \mathrm{CNY})\), while reducing the mean CVaR loss from \(27{,}208.41\) to \(11{,}191.13\) \((10^4\ \mathrm{CNY})\). The joint reserve reliability also increases from 0.7679 to 0.7931. These results indicate that second-stage appraisal projects provide not only additional value, but also an effective mechanism for compensating unfavorable first-stage outcomes.

The benefit of posterior updating can be observed by comparing posterior-informed exact recourse with fixed-probability exact recourse. The posterior-informed setting increases the mean expected NPV by \(5{,}462.11\) \((10^4\ \mathrm{CNY})\), corresponding to an improvement of 1.63\%, and reduces the mean CVaR loss by \(800.98\) \((10^4\ \mathrm{CNY})\), corresponding to a 6.68\% reduction. The mean posterior success probability of selected appraisal projects also rises from 0.8980 to 0.9139 (in Table.~\ref{Tab6}). Therefore, the logit-scale posterior update does not merely change probability labels; it changes the quality of the selected recourse projects by directing the second-stage optimization toward appraisal opportunities that are more consistent with the revealed geological information.
\begin{table*}[!t]
	\centering
	\caption{Out-of-sample performance under different recourse settings.}
	\label{Tab6}
	\begin{tabular}{lccccc}
		\hline
		Recourse setting & ENPV & CVaR loss & Success reliability & Reserve reliability & Selected posterior PoS \\
		\hline
		No recourse & 253732.31 & 27208.41 & 0.9981 & 0.7679 & -- \\
		Fixed-probability exact recourse & 334458.25 & 11992.11 & 0.9998 & 0.7928 & 0.8980 \\
		Posterior-informed exact recourse & 339920.36 & 11191.13 & 0.9999 & 0.7931 & 0.9139 \\
		Greedy recourse & 339599.93 & 11204.99 & 0.9999 & 0.7926 & 0.9128 \\
		\hline
	\end{tabular}
	\\
	\vspace{1mm}
	\footnotesize{Note: ENPV and CVaR loss are reported in \(10^4\ \mathrm{CNY}\). ``Selected posterior PoS'' denotes the mean posterior probability of success of selected second-stage appraisal projects.}
\end{table*}

\begin{figure}
	\centering
	\includegraphics[width=\linewidth]{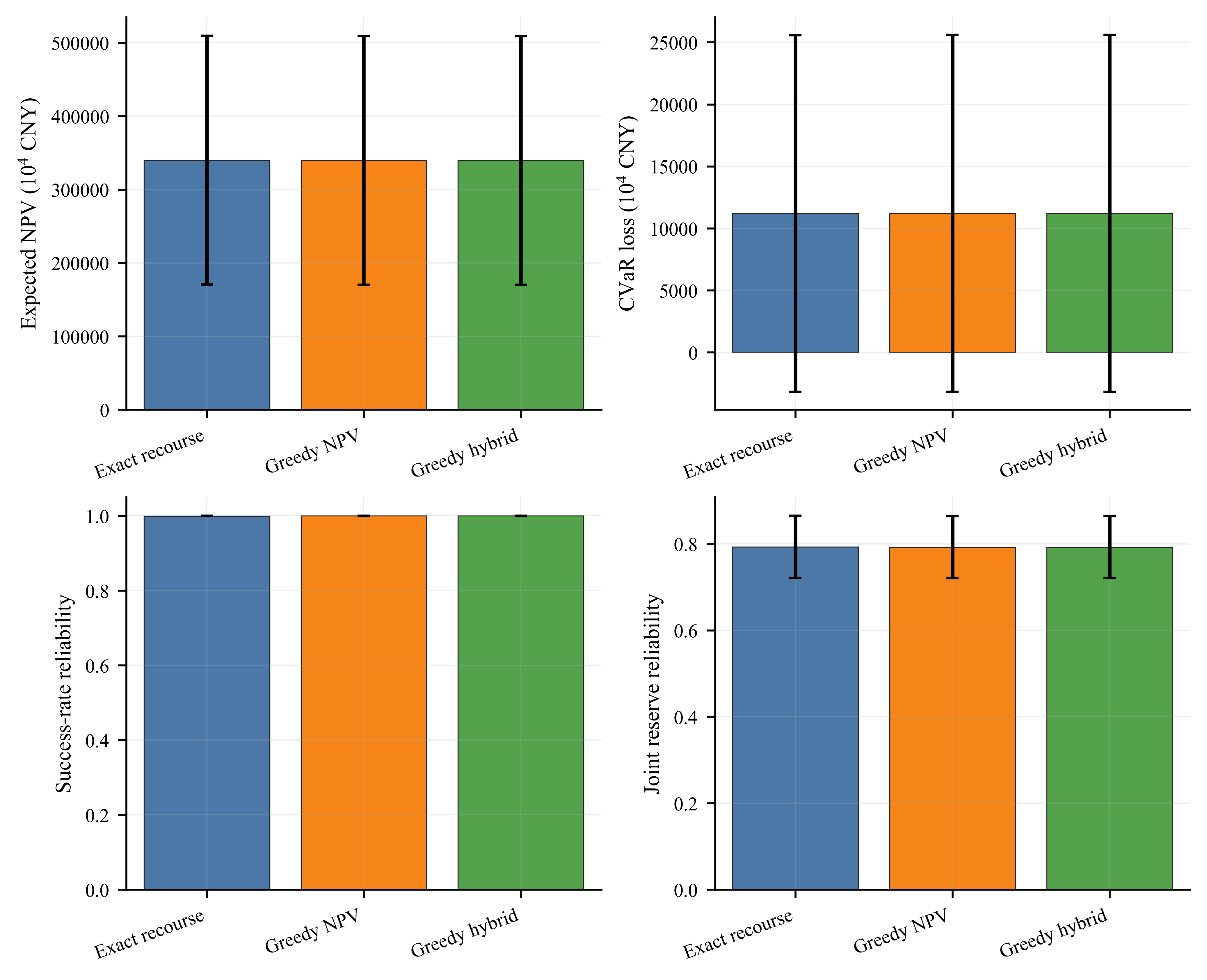}
	\caption{Out-of-sample comparison between exact recourse optimization and greedy recourse rules.}
	\label{fig7}
\end{figure}

The second experiment compares exact recourse optimization with greedy recourse rules. In the exact recourse setting, the second-stage decision is obtained by solving a scenario-wise constrained binary recourse problem, while the greedy settings select eligible appraisal projects according to local ranking rules. As shown in Fig.~\ref{fig7}, the exact recourse model produces a slightly higher mean expected NPV and a slightly lower mean CVaR loss than the greedy rules. Specifically, exact recourse improves the mean expected NPV by \(320.42\) \((10^4\ \mathrm{CNY})\) and reduces the mean CVaR loss by \(13.86\) \((10^4\ \mathrm{CNY})\) relative to greedy recourse.

Although the numerical difference between exact and greedy recourse is moderate in this case, the modeling implication is important. The exact recourse model evaluates second-stage appraisal projects jointly under remaining budget, well-count, eligibility, and project-category restrictions, whereas greedy rules make sequential local selections. Therefore, exact recourse provides a constraint-consistent and reproducible evaluation of each first-stage portfolio. This feature is essential for embedding recourse decisions within NSGA-II, because the quality of a first-stage portfolio should be assessed by its best feasible scenario-wise adjustment rather than by an ad hoc local selection rule.

\subsection{Out-of-sample robustness and sensitivity analysis}
This subsection evaluates the stability of the stochastic portfolios from two complementary perspectives. First, the portfolios obtained from the in-sample scenario bank are re-evaluated under an independent out-of-sample scenario bank. Second, sensitivity tests are conducted for the CVaR confidence level and the geological learning-strength coefficient. The purpose is to examine whether the proposed framework remains interpretable when the scenario sample, risk measurement, and posterior-learning intensity change.

\begin{figure}[!t]
	\centering
	\includegraphics[width=\linewidth]{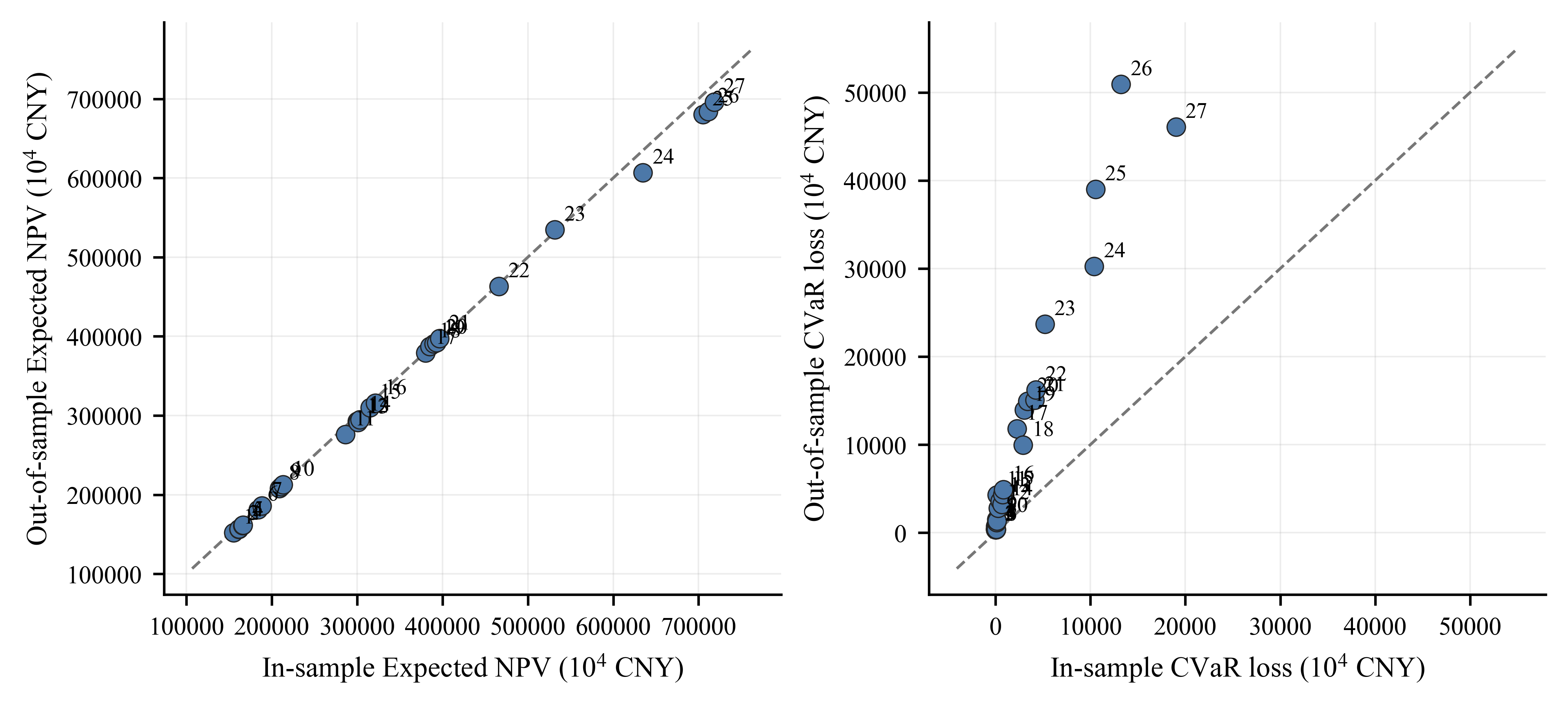}
	\caption{In-sample and out-of-sample performance of posterior-informed exact recourse portfolios.}
	\label{fig8}
\end{figure}

Fig.~\ref{fig8} compares the in-sample and out-of-sample performance of the posterior-informed exact recourse portfolios. The mean expected NPV decreases from \(346{,}395.45\) to \(339{,}920.36\) \((10^4\ \mathrm{CNY})\), corresponding to a reduction of approximately 1.87\%. This relatively small change indicates that the economic performance of the Pareto portfolios is not highly dependent on the in-sample scenario bank. At the same time, the mean CVaR loss increases from \(3{,}038.74\) to \(11{,}191.13\) \((10^4\ \mathrm{CNY})\). This increase is expected because the independent out-of-sample bank contains more scenarios and therefore provides a stricter test of tail losses.

The reliability indicators remain stable under out-of-sample validation. The mean success reliability increases from 0.9996 to 0.9999, and the mean joint reserve reliability changes from 0.7742 to 0.7931 (in Table.~\ref{Tab7}). These results suggest that the proposed model does not simply overfit economic outcomes in the in-sample scenarios; it also preserves geological and reserve-related reliability when evaluated on unseen uncertainty realizations.
\begin{table*}[!t]
	\centering
	\caption{Robustness and sensitivity results under out-of-sample evaluation.}
	\label{Tab7}
	\begin{tabular}{lccccc}
		\hline
		Experiment & Setting & ENPV & CVaR loss & Success reliability & Reserve reliability \\
		\hline
		Robustness & In-sample & 346395.45 & 3038.74 & 0.9996 & 0.7742 \\
		Robustness & Out-of-sample & 339920.36 & 11191.13 & 0.9999 & 0.7931 \\
		Risk sensitivity & \(\beta=0.80\) & 339920.36 & 5595.57 & 0.9999 & 0.7931 \\
		Risk sensitivity & \(\beta=0.90\) & 339920.36 & 11191.13 & 0.9999 & 0.7931 \\
		Risk sensitivity & \(\beta=0.95\) & 339920.36 & 22382.26 & 0.9999 & 0.7931 \\
		Learning sensitivity & \(\theta\)-scale \(=0.0\) & 334458.25 & 11992.11 & 0.9998 & 0.7928 \\
		Learning sensitivity & \(\theta\)-scale \(=1.0\) & 339920.36 & 11191.13 & 0.9999 & 0.7931 \\
		Learning sensitivity & \(\theta\)-scale \(=2.0\) & 344993.69 & 10439.47 & 0.9999 & 0.7933 \\
		\hline
	\end{tabular}
	\\
	\vspace{1mm}
	\footnotesize{Note: ENPV and CVaR loss are reported in \(10^4\ \mathrm{CNY}\). Risk sensitivity changes the CVaR confidence level for the same evaluated portfolios, whereas learning sensitivity changes the posterior-learning strength used in recourse evaluation.}
\end{table*}

\begin{figure}[!t]
	\centering
	\includegraphics[width=\linewidth]{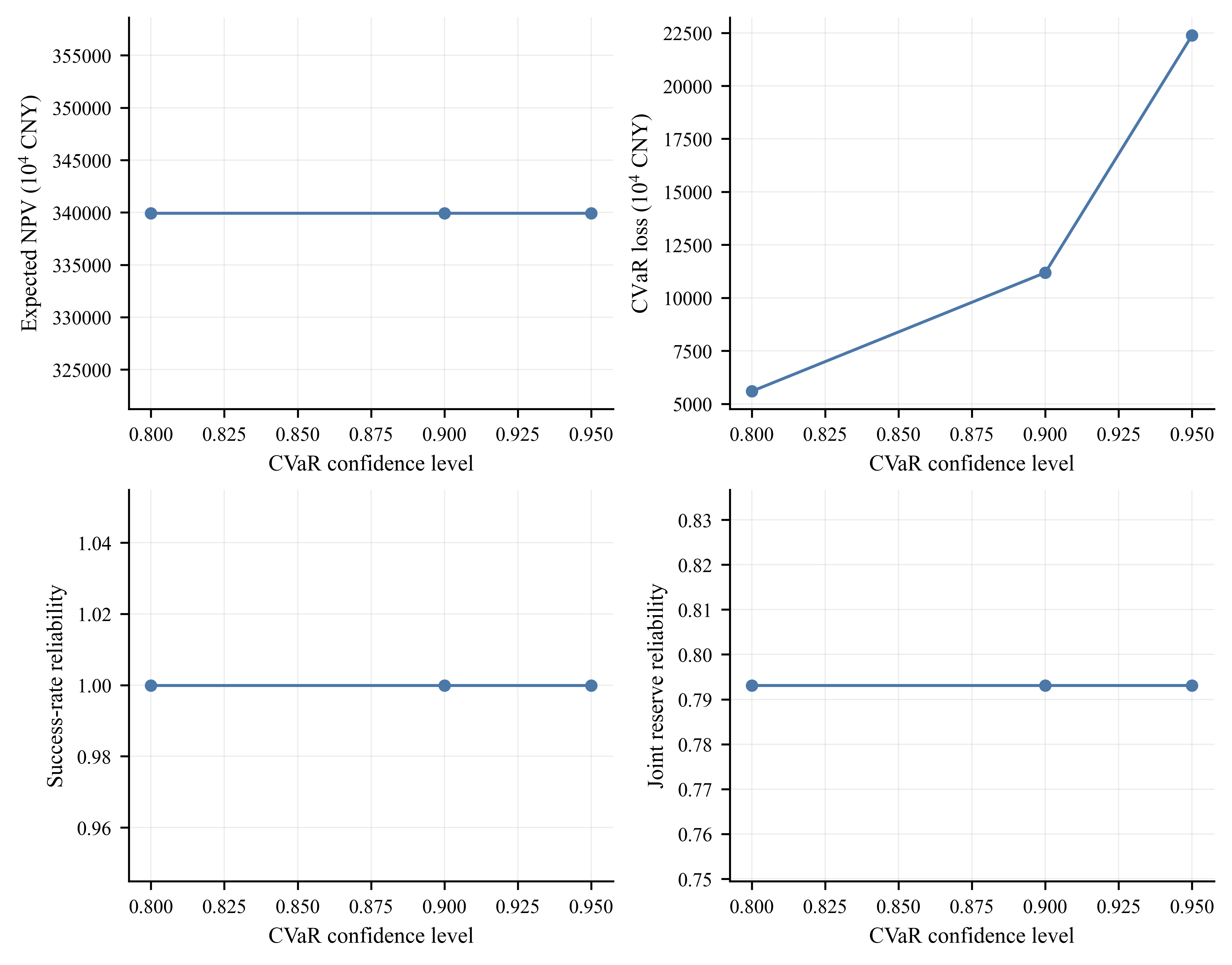}
	\caption{Sensitivity of portfolio performance to the CVaR confidence level.}
	\label{fig9}
\end{figure}

Fig.~\ref{fig9} reports the sensitivity of the CVaR confidence level. When \(\beta\) increases from 0.80 to 0.95, the mean CVaR loss increases from \(5{,}595.57\) to \(22{,}382.26\) \((10^4\ \mathrm{CNY})\). The expected NPV and reliability indicators remain unchanged because this experiment re-evaluates the same portfolio set under different CVaR confidence levels rather than re-optimizing the portfolio under each \(\beta\). Therefore, the result should be interpreted as a diagnostic of tail-risk exposure: a higher confidence level reveals more severe downside losses in the extreme tail of the NPV distribution.

\begin{figure}[!t]
	\centering
	\includegraphics[width=\linewidth]{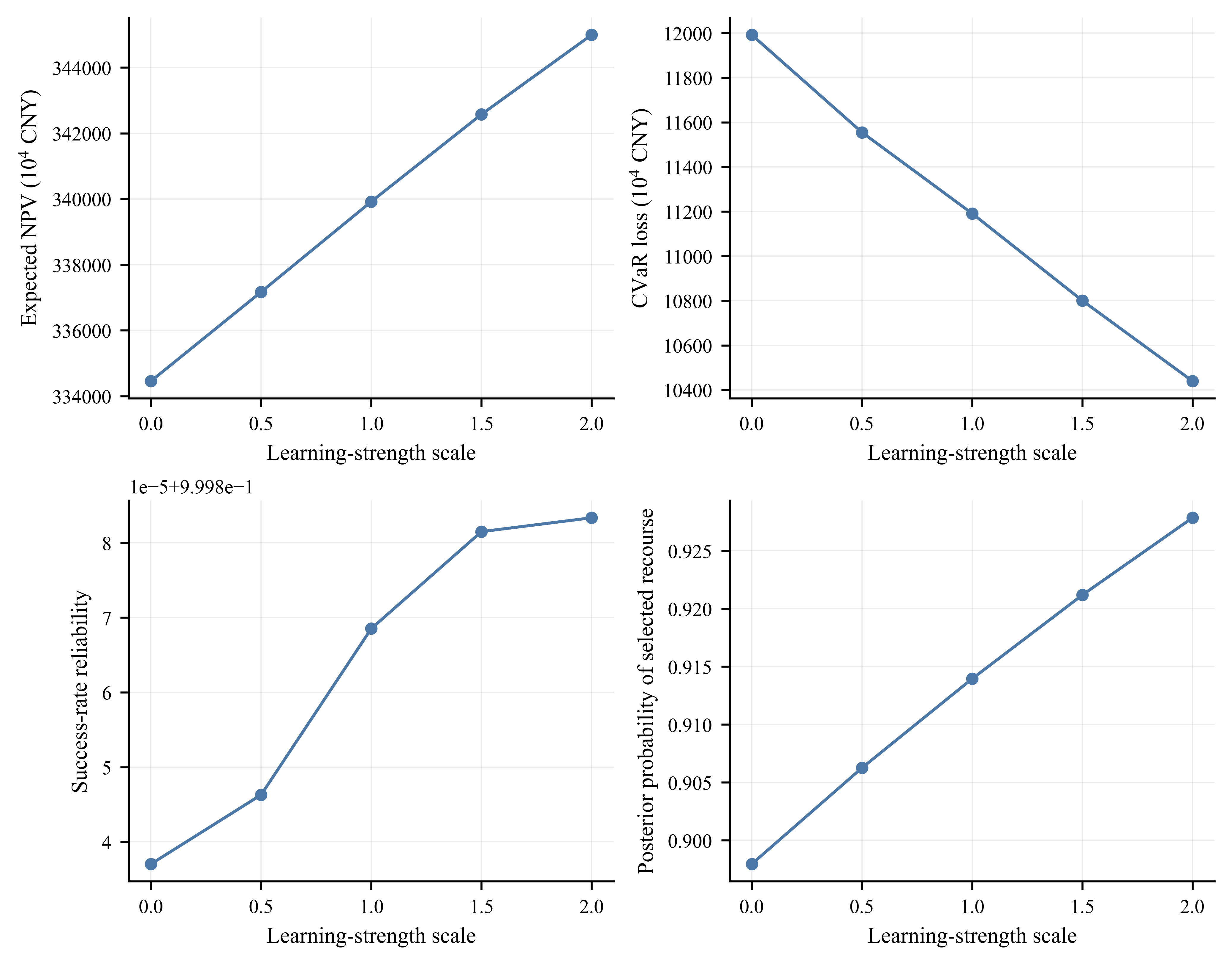}
	\caption{Sensitivity of portfolio performance to geological learning strength.}
	\label{fig10}
\end{figure}

Fig.~\ref{fig10} shows the effect of geological learning strength. As the learning-strength scale increases from 0.0 to 2.0, the mean expected NPV increases from \(334{,}458.25\) to \(344{,}993.69\) \((10^4\ \mathrm{CNY})\), while the mean CVaR loss decreases from \(11{,}992.11\) to \(10{,}439.47\) \((10^4\ \mathrm{CNY})\). The mean posterior probability of selected appraisal projects also rises from 0.8980 to 0.9278. These monotonic changes support the role of posterior-informed learning: stronger geological information transfer improves the selection quality of second-stage appraisal projects and reduces downside exposure without undermining reliability.

\begin{figure}[!t]
	\centering
	\includegraphics[width=\linewidth]{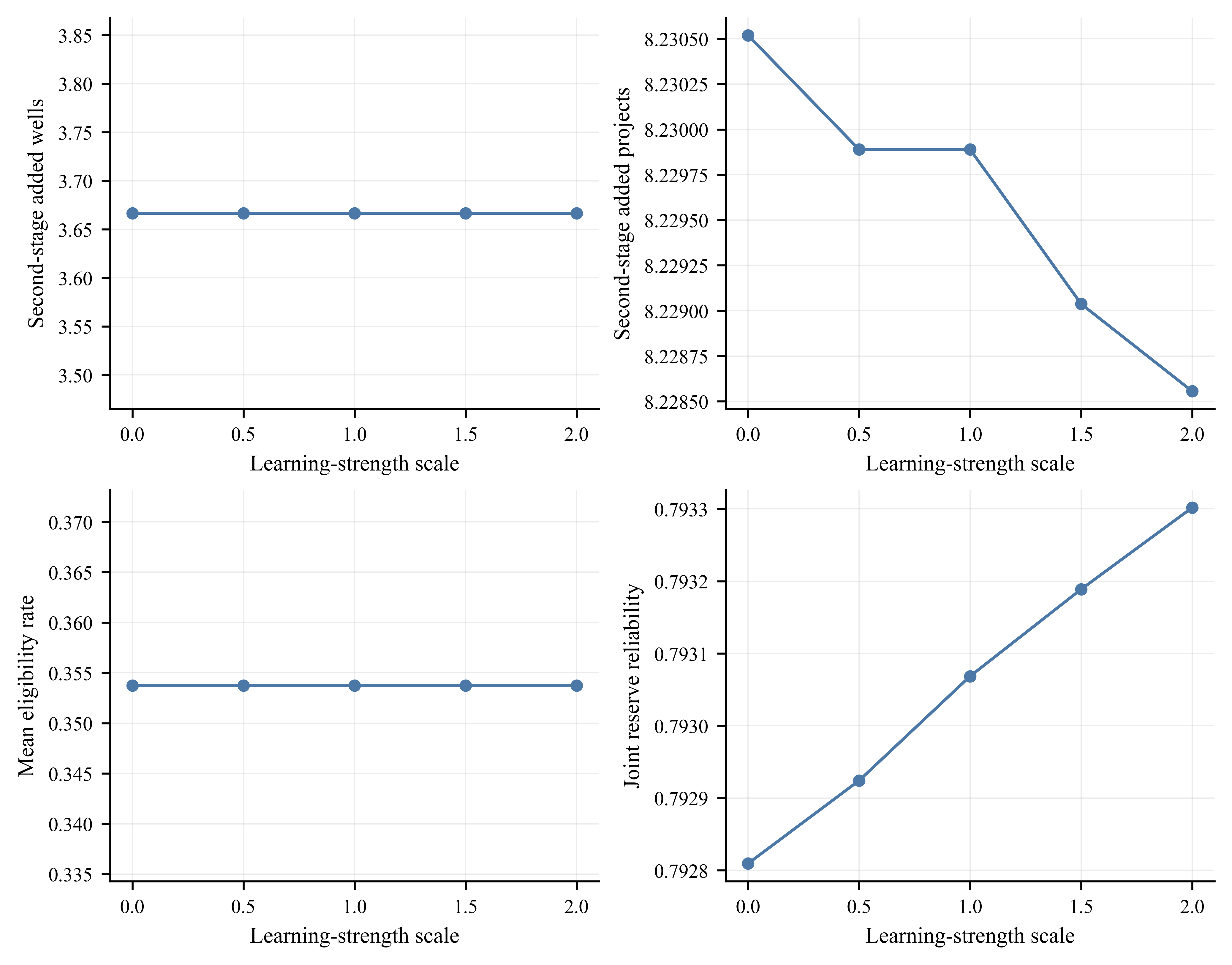}
	\caption{Second-stage recourse behavior under different geological learning-strength settings.}
	\label{fig11}
\end{figure}

The recourse-behavior indicators in Fig.~\ref{fig11} provide an additional explanation. The average number of second-stage wells remains nearly unchanged across learning-strength settings, while the selected posterior probability improves. This means that the performance gain is not mainly caused by drilling more wells or expanding recourse scale; it is driven by better information-guided selection among feasible appraisal projects. This distinction is important for exploration management because it links geological learning directly to portfolio quality rather than to a simple increase in investment intensity.

\subsection{Numerical reliability of scenario sampling and evolutionary search}
The preceding subsections evaluate the economic value, recourse mechanism, robustness, and sensitivity of the proposed framework. This subsection further examines the numerical reliability of the computational procedure from two perspectives: the stability of SAA-based scenario evaluation and the repeatability of the NSGA-II search under different random seeds. These tests are necessary because both scenario sampling and evolutionary optimization introduce numerical variability into the obtained Pareto frontier.

\begin{figure}[!t]
	\centering
	\includegraphics[width=\linewidth]{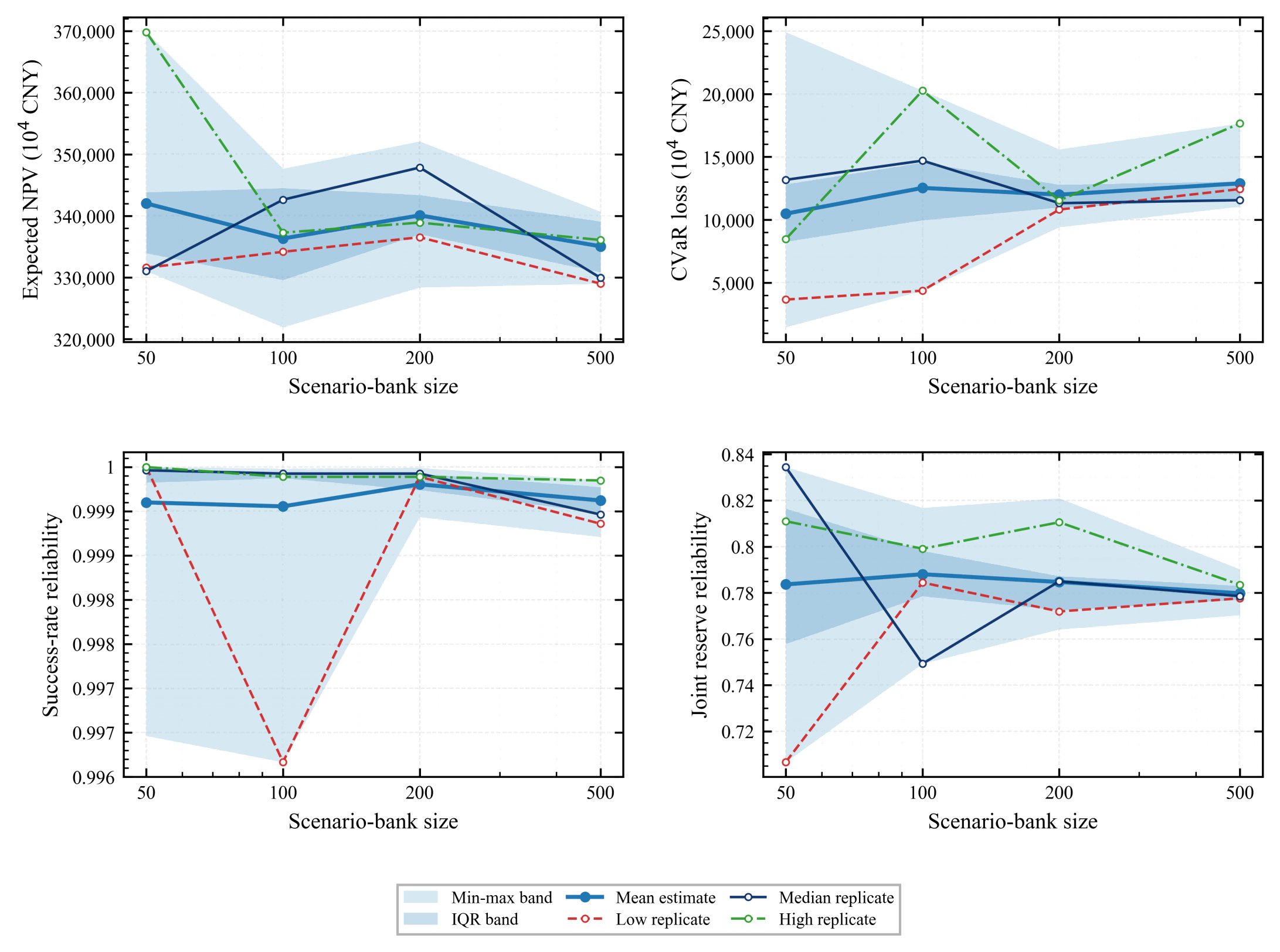}
	\caption{Scenario stability of SAA-based performance estimates under different scenario-bank sizes.}
	\label{fig12}
\end{figure}

The SAA stability test evaluates the posterior-informed exact recourse portfolios under scenario banks with \(S=50\), \(100\), \(200\), and \(500\). For each scenario size, 10 independent scenario banks are generated, while the number of second-stage sub-scenarios is fixed at \(K=20\). Fig.~\ref{fig12} reports the min--max band, interquartile range (IQR) band, mean estimate, and representative low-, median-, and high-replicate trajectories for each performance indicator. The expected NPV estimate becomes more concentrated as the scenario size increases: the coefficient of variation of mean expected NPV decreases from 3.26\% at \(S=50\) to 1.32\% at \(S=500\). The joint reserve-reliability estimate also becomes more stable, with its standard deviation decreasing from 0.0418 to 0.0064. These results indicate that increasing the scenario-bank size reduces sampling-induced variability in both economic and reserve-reliability evaluation.
\begin{table*}[!t]
	\centering
	\caption{Numerical stability of SAA evaluation and NSGA-II repeated search.}
	\label{Tab8}
	\begin{tabular}{lcccc}
		\hline
		Test & Indicator & Minimum & Mean & Maximum \\
		\hline
		SAA, \(S=50\) & CV of ENPV & -- & 3.26\% & -- \\
		SAA, \(S=500\) & CV of ENPV & -- & 1.32\% & -- \\
		SAA, \(S=50\) & CV of CVaR & -- & 61.15\% & -- \\
		SAA, \(S=500\) & CV of CVaR & -- & 14.61\% & -- \\
		NSGA-II repeated runs & Feasible rate & 1.0000 & 1.0000 & 1.0000 \\
		NSGA-II repeated runs & Front size & 9 & 15.3 & 23 \\
		NSGA-II repeated runs & Best ENPV & 718529.11 & 734001.60 & 737836.26 \\
		NSGA-II repeated runs & Minimum CVaR loss & 0.00 & 4.75 & 17.78 \\
		NSGA-II repeated runs & Hypervolume & \(1.2543\times 10^9\) & \(1.8328\times 10^9\) & \(2.1116\times 10^9\) \\
		\hline
	\end{tabular}
	\\
	\vspace{1mm}
	\footnotesize{Note: ENPV and CVaR loss are reported in \(10^4\ \mathrm{CNY}\). The SAA coefficients of variation are calculated across 10 independent scenario banks for each scenario size. The NSGA-II statistics are calculated across 10 independent runs with different random seeds.}
\end{table*}

The CVaR estimate is more sensitive to scenario sampling than the expected NPV estimate. Its coefficient of variation decreases from 61.15\% at \(S=50\) to 14.61\% at \(S=500\) (Table.~\ref{Tab8}) but remains substantially higher than that of expected NPV. This result is consistent with the tail-risk nature of CVaR, since rare adverse outcomes require larger scenario banks to be adequately represented. Therefore, the SAA experiment confirms that larger scenario samples improve numerical stability, while also explaining why independent out-of-sample validation is necessary for assessing downside-risk robustness.

To assess the repeatability of the evolutionary search, NSGA-II is independently executed 10 times with different random seeds. Each run uses the same population size, number of generations, and evaluation scenario bank. All repeated runs return feasible portfolios, and the feasible rate remains 1.0000. The number of non-dominated portfolios varies from 9 to 23, with an average of 15.3, indicating that the density and coverage of the final frontier are affected by the stochastic nature of evolutionary search. However, the best expected NPV remains stable across runs, ranging from \(718{,}529.11\) to \(737{,}836.26\) \((10^4\ \mathrm{CNY})\), and the minimum CVaR loss remains low, with a mean of \(4.75\) \((10^4\ \mathrm{CNY})\).

\begin{figure}[!t]
	\centering
	\includegraphics[width=\linewidth]{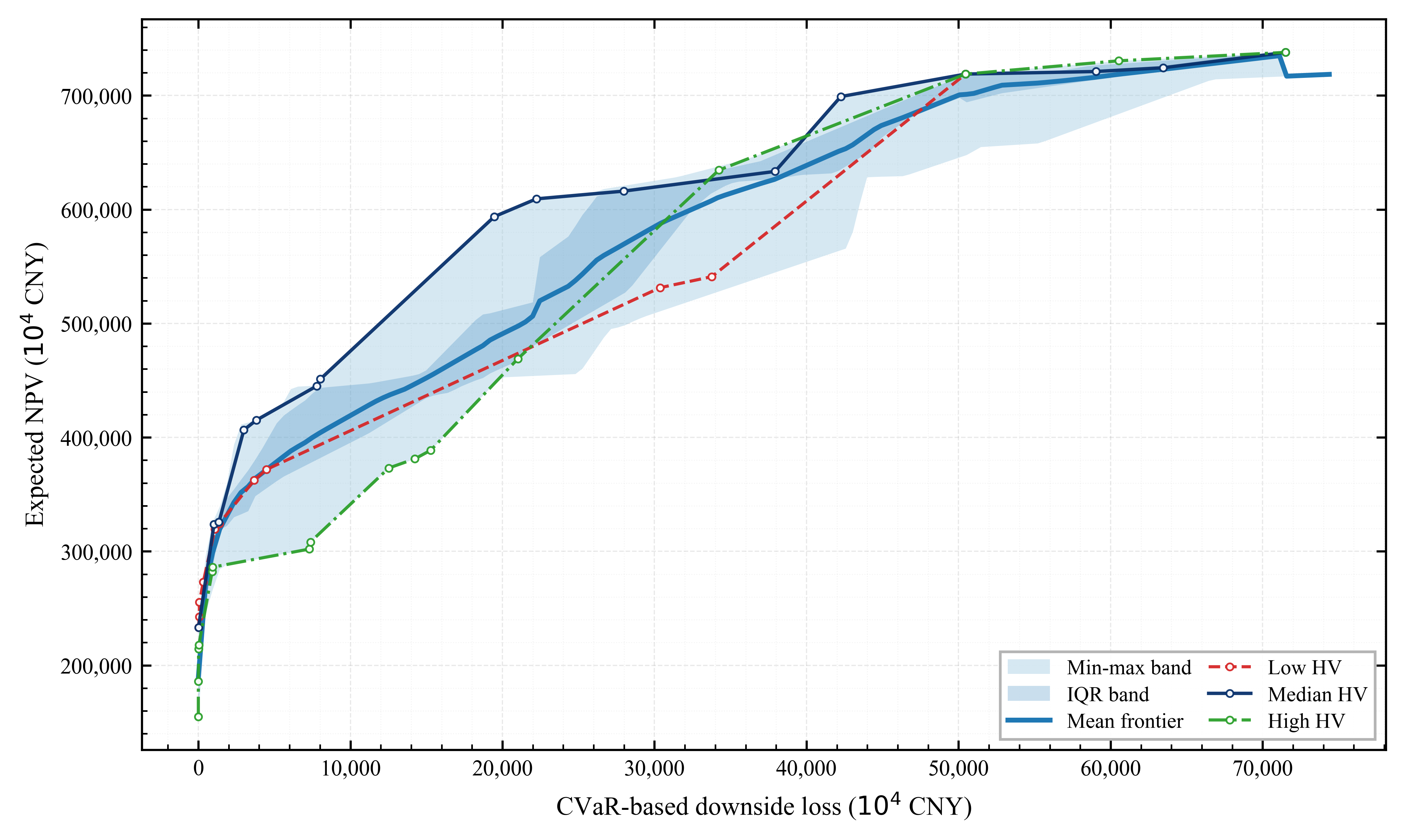}
	\caption{Variability of Pareto fronts obtained from repeated NSGA-II runs under different random seeds.}
	\label{fig13}
\end{figure}

Fig.~\ref{fig13} visualizes the variability of the Pareto fronts obtained from repeated NSGA-II runs. The min--max band describes the full range of interpolated frontiers, and the IQR band describes the central dispersion. The mean frontier remains within a stable trade-off region, while the low-, median-, and high-hypervolume runs show that the main source of variability lies in the coverage and density of the frontier rather than in the feasibility of the solutions. This result supports the numerical credibility of the proposed solution procedure: NSGA-II consistently identifies feasible and high-quality Pareto portfolios, although the exact distribution of non-dominated points may vary across random seeds.

Overall, the numerical reliability experiments indicate that the proposed framework is not dependent on a single scenario sample or a single evolutionary run. Increasing the scenario-bank size improves the stability of SAA estimates, especially for reserve reliability and expected NPV, while repeated NSGA-II runs confirm the feasibility and quality of the generated Pareto portfolios. These findings strengthen the computational basis of the proposed two-stage stochastic multi-objective optimization framework.

\section{Conclusion}\label{sec6}
This paper develops a posterior-informed two-stage stochastic multi-objective optimization framework for oil and gas exploration portfolio selection. The proposed model jointly considers geological success uncertainty, reserve uncertainty, and economic value uncertainty. First-stage decisions determine the initial exploration portfolio before drilling outcomes are observed, while second-stage appraisal projects are selected after scenario-specific geological information is revealed. A logit-scale posterior updating mechanism is used to link first-stage drilling outcomes with recourse-project success probabilities, and an explicit binary recourse optimization problem is solved under each scenario to evaluate feasible second-stage adjustments. The resulting first-stage portfolios are searched by NSGA-II and evaluated by expected NPV, CVaR-based downside risk, geological success reliability, and reserve reliability.

The computational results show that the proposed framework generates a meaningful return--risk frontier. The model obtains 27 feasible non-dominated portfolios, with out-of-sample expected NPV ranging from \(151{,}967.48\) to \(695{,}884.35\) \((10^4\ \mathrm{CNY})\) and CVaR-based downside loss ranging from \(337.62\) to \(50{,}934.16\) \((10^4\ \mathrm{CNY})\). The success reliability remains high across the frontier, whereas the joint reserve reliability varies substantially, indicating that geological success and reserve deliverability should be treated as distinct decision criteria.

The comparison with a deterministic mean-value benchmark further demonstrates the necessity of stochastic optimization. Although the deterministic benchmark achieves a high out-of-sample expected NPV, it produces a larger CVaR loss and violates feasibility requirements under stochastic evaluation. This result suggests that deterministic expected-value optimization may overstate portfolio attractiveness when geological, reserve, and economic uncertainties are ignored.

The recourse experiments confirm the value of posterior-informed decision adjustment. Compared with the no-recourse setting, posterior-informed exact recourse increases the mean expected NPV by 33.97\% and reduces the mean CVaR loss by 58.87\%. Compared with fixed-probability recourse, posterior-informed recourse increases the mean expected NPV by 1.63\% and reduces the mean CVaR loss by 6.68\%. These findings show that revealed geological information can improve second-stage project selection rather than merely adding flexibility. The exact recourse model also provides a constraint-consistent evaluation of second-stage decisions, even though its numerical advantage over greedy recourse is moderate in the present case.

Out-of-sample validation and sensitivity analysis support the robustness of the proposed framework. The mean expected NPV decreases by only 1.87\% from in-sample to out-of-sample evaluation, while reliability indicators remain stable. The CVaR estimate is more sensitive to scenario sampling, which confirms the importance of explicitly modeling tail risk. In addition, increasing geological learning strength improves the selected posterior probability, raises expected NPV, and reduces CVaR loss without increasing the average number of second-stage wells. This indicates that the benefit of geological learning is mainly driven by better information-guided recourse selection rather than by simply expanding drilling intensity.

Overall, the proposed framework provides a structured decision-support tool for exploration portfolio planning under uncertainty. It integrates stochastic scenario modeling, posterior geological learning, explicit recourse optimization, and multi-objective portfolio search into a unified framework. Future work may extend the model by calibrating uncertainty distributions with field data, incorporating spatial geological dependence among prospects, and considering multi-period exploration strategies with sequential information acquisition.

\section*{Credit authorship contribution statement}
\textbf{Junyi Cui:} Writing – original draft, Investigation, Resources, Methodology, Formal analysis. \textbf{Junyi Cui}: Formal analysis, Software, Data curation, Writing – review \& editing, Validation. \textbf{Junyi Cui}: Validation, Formal analysis, Visualization, Supervision. \textbf{Junyi Cui}: Formal analysis, Validation, Supervision. \textbf{Junyi Cui}: Validation, Writing – review \& editing, Supervision.
\section*{Declaration of competing interest}

The authors declare that they have no known competing financial interests or personal relationships that could have appeared to influence the work reported in this paper.

\section*{Acknowledgments}
This work was supported by Sichuan Provincial Science and Technology Project (No.2025NSFTD0016) and the Open Bidding for Selecting the Best Candidates Project of Southwest Petroleum University (2024CXJB11).


\section*{Data availability}
The data that has been used is confidential.

\bibliographystyle{elsarticle-num}
\bibliography{ces-refs}

@article{bib1,
title = {Petroleum exploration portfolios generated with different optimization approaches: Lessons for decision-makers},
journal = {Journal of Petroleum Science and Engineering},
volume = {214},
pages = {110459},
year = {2022},
issn = {0920-4105},
doi = {https://doi.org/10.1016/j.petrol.2022.110459},
author = {William C.P. LaCosta and Alexei V. Milkov}
}

@article{bib2,
title = {Risk analysis applied to petroleum exploration and production: an overview},
journal = {Journal of Petroleum Science and Engineering},
volume = {44},
number = {1},
pages = {1-9},
year = {2004},
note = {Risk Analysis Applied to Petroleum Exploration and Production},
issn = {0920-4105},
doi = {https://doi.org/10.1016/j.petrol.2004.02.001},
author = {S.B. Suslick and D.J. Schiozer}
}

@article{bib3,
title = {Risk tables for less biased and more consistent estimation of probability of geological success (PoS) for segments with conventional oil and gas prospective resources},
journal = {Earth-Science Reviews},
volume = {150},
pages = {453-476},
year = {2015},
issn = {0012-8252},
doi = {https://doi.org/10.1016/j.earscirev.2015.08.006},
author = {Alexei V. Milkov}
}

@article{bib4,
title = {Combining decision analysis and portfolio management to improve project selection in the exploration and production firm},
journal = {Journal of Petroleum Science and Engineering},
volume = {44},
number = {1},
pages = {55-65},
year = {2004},
note = {Risk Analysis Applied to Petroleum Exploration and Production},
issn = {0920-4105},
doi = {https://doi.org/10.1016/j.petrol.2004.02.005},
author = {Michael R. Walls}
}

@article{bib5,
title = {Expected value, downside risk and upside potential as decision criteria in production strategy selection for petroleum field development},
journal = {Journal of Petroleum Science and Engineering},
volume = {157},
pages = {81-93},
year = {2017},
issn = {0920-4105},
doi = {https://doi.org/10.1016/j.petrol.2017.07.002},
author = {Susana M.G. Santos and Vinicius E. Botechia and Denis J. Schiozer and Ana T.F.S. Gaspar},
}

@article{bib6,
title = {An exposition of least square Monte Carlo approach for real options valuation},
journal = {Geoenergy Science and Engineering},
volume = {222},
pages = {111230},
year = {2023},
issn = {2949-8910},
doi = {https://doi.org/10.1016/j.petrol.2022.111230},
author = {Rouholah Ahmadi and Reidar Brumer Bratvold}
}

@article{bib7,
title = {A multi-objective decision-making framework for the choice between mutually exclusive alternatives under uncertainty: Assessing the competitiveness of offshore wind for a gas field electrification on the NCS},
journal = {Energy Economics},
volume = {141},
pages = {108032},
year = {2025},
issn = {0140-9883},
doi = {https://doi.org/10.1016/j.eneco.2024.108032},
author = {Daniel Aghajani and Reidar B. Bratvold and Verena Hagspiel and Olga Noshchenko and Vincent K.G. Toutain},
}

@Article{bib8,
AUTHOR = {Yong, Chi and Tong, Mu and Yang, Zhongyi and Zhou, Jixian},
TITLE = {Conventional Natural Gas Project Investment and Decision Making under Multiple Uncertainties},
JOURNAL = {Energies},
VOLUME = {16},
YEAR = {2023},
NUMBER = {5},
ARTICLE-NUMBER = {2342},
ISSN = {1996-1073},
DOI = {10.3390/en16052342}
}

@article{bib9,
  title={How to optimize the investment decision of shale gas multi-objective green development for energy sustainable development?},
  author={Liu, Haomin and Zhang, Tao},
  journal={Journal of Petroleum Exploration and Production Technology},
  volume={15},
  number={12},
  pages={192},
  year={2025},
  publisher={Springer},
  doi={https://doi.org/10.1007/s13202-025-02084-7}
}

@article{bib10,
title = {A review of the operations literature on real options in energy},
journal = {European Journal of Operational Research},
volume = {309},
number = {2},
pages = {469-487},
year = {2023},
issn = {0377-2217},
doi = {https://doi.org/10.1016/j.ejor.2022.09.014},
author = {Selvaprabu Nadarajah and Nicola Secomandi}
}

@article{bib11,
title = {Environmental and economic multi-objective real options analysis: Electrification choices for field development investment planning},
journal = {Energy},
volume = {295},
pages = {131053},
year = {2024},
issn = {0360-5442},
doi = {https://doi.org/10.1016/j.energy.2024.131053},
author = {Olga Noshchenko and Verena Hagspiel}
}

@article{bib12,
title = {Valuing oil reserve volumes under price uncertainty},
journal = {Journal of Economics and Business},
volume = {137},
pages = {106277},
year = {2025},
note = {Recent Developments and Challenges in Real Options},
issn = {0148-6195},
doi = {https://doi.org/10.1016/j.jeconbus.2025.106277},
author = {Marco Antonio Guimarães Dias and Roberto Evelim Penha Borges}
}

@article{bib13,
  title={A multicriteria decision model for selecting a portfolio of oil and gas exploration projects},
  author={Lopes, Yuri Gama and Almeida, Adiel Teixeira de},
  journal={Pesquisa Operacional},
  volume={33},
  number={3},
  pages={417--441},
  year={2013},
  publisher={SciELO Brasil},
  doi={https://doi.org/10.1590/S0101-74382013005000011}
}

@article{bib14,
title = {Assessment of synergies for selecting a project portfolio in the petroleum industry based on a multi-attribute utility function},
journal = {Journal of Petroleum Science and Engineering},
volume = {126},
pages = {131-140},
year = {2015},
issn = {0920-4105},
doi = {https://doi.org/10.1016/j.petrol.2014.12.012},
author = {Yuri Gama Lopes and Adiel Teixeira {de Almeida}},
}

@article{bib15,
  title={Stochastic programming approach for the planning of offshore oil or gas field infrastructure under decision-dependent uncertainty},
  author={Tarhan, Bora and Grossmann, Ignacio E and Goel, Vikas},
  journal={Industrial \& Engineering Chemistry Research},
  volume={48},
  number={6},
  pages={3078--3097},
  year={2009},
  publisher={ACS Publications}
}

@article{bib16,
title = {Multistage stochastic programming approach for offshore oilfield infrastructure planning under production sharing agreements and endogenous uncertainties},
journal = {Journal of Petroleum Science and Engineering},
volume = {124},
pages = {180-197},
year = {2014},
issn = {0920-4105},
doi = {https://doi.org/10.1016/j.petrol.2014.10.006},
author = {Vijay Gupta and Ignacio E. Grossmann}
}

@article{bib17,
    author = {Martinelli, Gabriele  and Eidsvik, Jo  and Hokstad, Ketil  and Hauge, Ragnar },
    title = {Strategies for Petroleum Exploration on the Basis of Bayesian Networks: A Case Study},
    journal = {SPE Journal},
    volume = {19},
    number = {04},
    pages = {564-575},
    year = {2013},
    month = {08},
    doi = {10.2118/159722-PA}
}

@article{bib18,
title = {Dynamic decision making for graphical models applied to oil exploration},
journal = {European Journal of Operational Research},
volume = {230},
number = {3},
pages = {688-702},
year = {2013},
issn = {0377-2217},
doi = {https://doi.org/10.1016/j.ejor.2013.04.057},
author = {Gabriele Martinelli and Jo Eidsvik and Ragnar Hauge}
}

@article{bib19,
title = {Multistage stochastic programming approach for offshore oilfield infrastructure planning under production sharing agreements and endogenous uncertainties},
journal = {Journal of Petroleum Science and Engineering},
volume = {124},
pages = {180-197},
year = {2014},
issn = {0920-4105},
doi = {https://doi.org/10.1016/j.petrol.2014.10.006},
author = {Vijay Gupta and Ignacio E. Grossmann}
}

@article{bib20,
title = {Data-driven project portfolio selection: Decision-dependent stochastic programming formulations with reliability and time to market requirements},
journal = {Computers \& Operations Research},
volume = {143},
pages = {105737},
year = {2022},
issn = {0305-0548},
doi = {https://doi.org/10.1016/j.cor.2022.105737},
author = {Janne Kettunen and Miguel A. Lejeune}
}

@article{bib21,
author = {Robert E. Kass and Adrian E. Raftery},
title = {Bayes Factors},
journal = {Journal of the American Statistical Association},
volume = {90},
number = {430},
pages = {773--795},
year = {1995},
publisher = {Taylor \& Francis},
doi = {10.1080/01621459.1995.10476572}
}

@article{bib22,
author = {Artzner, Philippe and Delbaen, Freddy and Eber, Jean-Marc and Heath, David},
title = {Coherent Measures of Risk},
journal = {Mathematical Finance},
volume = {9},
number = {3},
pages = {203-228},
doi = {https://doi.org/10.1111/1467-9965.00068},
year = {1999}
}

@article{bib23,
  title={Optimization of conditional value-at-risk},
  author={Rockafellar, R Tyrrell and Uryasev, Stanislav and others},
  journal={Journal of risk},
  volume={2},
  pages={21--42},
  year={2000},
  doi={10.21314/JOR.2000.038}
}

@article{bib24,
  title={Portfolio optimization with conditional value-at-risk objective and constraints},
  author={Krokhmal, Pavlo and Palmquist, Jonas and Uryasev, Stanislav},
  journal={Journal of risk},
  volume={4},
  pages={43--68},
  year={2002}
}

@article{bib25,
author = {Kleywegt, Anton J. and Shapiro, Alexander and Homem-de-Mello, Tito},
title = {The Sample Average Approximation Method for Stochastic Discrete Optimization},
journal = {SIAM Journal on Optimization},
volume = {12},
number = {2},
pages = {479-502},
year = {2002},
doi = {10.1137/S1052623499363220}
}

@ARTICLE{bib26,
  author={Deb, K. and Pratap, A. and Agarwal, S. and Meyarivan, T.},
  journal={IEEE Transactions on Evolutionary Computation}, 
  title={A fast and elitist multiobjective genetic algorithm: NSGA-II}, 
  year={2002},
  volume={6},
  number={2},
  pages={182-197},
  doi={10.1109/4235.996017}
  }

@article{bib27,
author = {Robert E. Kass and Adrian E. Raftery},
title = {Bayes Factors},
journal = {Journal of the American Statistical Association},
volume = {90},
number = {430},
pages = {773--795},
year = {1995},
publisher = {Taylor \& Francis},
doi = {10.1080/01621459.1995.10476572}
}

@article{bib28,
  author  = {Agterberg, Frederik P. and Cheng, Qiuming},
  title   = {Conditional Independence Test for Weights-of-Evidence Modeling},
  journal = {Natural Resources Research},
  volume  = {11},
  number  = {4},
  pages   = {249--255},
  year    = {2002},
  doi     = {10.1023/A:1021193827501}
}

@article{bib29,
  author  = {Monigle, P. W. and Hedayati, T. S. and Goulding, F. J.},
  title   = {Integrated and improved direct hydrocarbon indicators: A step forward in petroleum risk discrimination},
  journal = {AAPG Bulletin},
  volume  = {109},
  number  = {5},
  pages   = {617--636},
  year    = {2025},
  doi     = {10.1306/04042524030}
}

@article{bib30,
title = {Conditional value-at-risk for general loss distributions},
journal = {Journal of Banking \& Finance},
volume = {26},
number = {7},
pages = {1443-1471},
year = {2002},
issn = {0378-4266},
doi = {https://doi.org/10.1016/S0378-4266(02)00271-6},
author = {R.Tyrrell Rockafellar and Stanislav Uryasev}
}
\end{document}